\documentclass[12pt]{amsart}

\usepackage{graphicx}
\usepackage{epsfig}
\usepackage{pstricks}

\newtheorem{thm}{Theorem}[section]

\newtheorem{defn}[thm]{Definition}
\newtheorem{lemma}[thm]{Lemma}

\newtheorem{cor}[thm]{Corollary}

\newtheorem{example}[thm]{Example}

\usepackage{amsmath}
\usepackage{amsxtra}
\usepackage{amscd}
\usepackage{amsthm}
\usepackage{amsfonts}
\usepackage{amssymb}
\usepackage{eucal}
\newcommand{\cQ}{\mathcal Q}
\newcommand{\cT}{\mathcal T}

\newcommand{\cC}{\mathcal C}
\newcommand{\g}{{\mathfrak{g}}}

\newcommand{\cP}{\mathcal P}

\newcommand{\Z}{{\mathbb Z}}

\newcommand{\ba}{{\mathbf a}}
\newcommand{\bb}{{\mathbf b}}

\newcommand{\bm}{{\mathbf m}}

\newcommand{\bx}{{\mathbf x}}
\newcommand{\by}{{\mathbf y}}

\newcommand{\al}{{\alpha}}

\numberwithin{equation}{section}

\begin{document}
\title{The solution of the $A_r$ T-system for arbitrary boundary}
\author{Philippe Di Francesco} 
\address{Department of Mathematics, University of Michigan,
530 Church Street, Ann Arbor, MI 48190, USA
and Institut de Physique Th\'eorique du Commissariat \`a l'Energie Atomique, 
Unit\'e de Recherche associ\'ee du CNRS,
CEA Saclay/IPhT/Bat 774, F-91191 Gif sur Yvette Cedex, 
FRANCE. e-mail: philippe.di-francesco@cea.fr}

\begin{abstract}
We present an explicit solution of the $A_r$ $T$-system for arbitrary boundary conditions.
For each boundary,
this is done by constructing a network, i.e. a graph with positively weighted edges,
and the solution is expressed as the partition function for a family of
non-intersecting paths on the network. This proves in particular the positive Laurent
property, namely that the solutions are all Laurent polynomials of the initial data with
non-negative integer coefficients.
\end{abstract}

\maketitle
\date{\today}
\tableofcontents

\section{Introduction}

In this paper we study the solutions of the $A_r$ $T$-system, namely
the following coupled system of recursion relations for $\al,j,k\in \Z$:
\begin{equation}\label{Tsys}
T_{\al,j,k+1}T_{\al,j,k-1}=T_{\al,j+1,k}T_{\al,j-1,k}+T_{\al+1,j,k}T_{\al-1,j,k}
\end{equation}
for $\al\in  I_r=\{1,2,...,r\}$, and
subject to the boundary conditions
\begin{equation}\label{bounT}
T_{0,j,k}=T_{r+1,j,k}=1 \qquad (j,k\in\Z)
\end{equation}

This system arose in many different contexts.
The system \eqref{Tsys} and its generalizations were introduced as the set of relations
satisfied by the eigenvalues 
of the fused transfer matrices of generalized quantum spin chains
based on any simply-laced Lie algebra $\g$  \cite{BR} \cite{KNS}; 
in this paper we restrict ourselves to the case $\g=sl_{r+1}$, 
but we believe our constructions can be adapted to other $\g$'s as well.

With the additional condition
that $T_{\al,0,k}=1$, $k\in\Z$ and the restriction to $j\in\Z_+$,
the solutions of (\ref{Tsys}-\ref{bounT}) were
also interpreted as the q-characters of some representations of the
affine Lie algebra $U_q(\widehat{sl}_{r+1})$,
the so-called Kirillov-Reshetikhin modules  \cite{FR},
indexed by $\al\in I_r=\{1,2,...,r\}$ and $j\in\Z_+$, 
while $k$ stands for a discrete spectral parameter \cite{Nakajima}.

The same equations appeared in the context of enumeration of domino tilings of plane domains
\cite{SPY}, and was studied in its own right under the name of octahedron equation \cite{KT} \cite{AH}. 
As noted by many authors, this equation may also be viewed
as a particular case of Pl\"ucker relations when all $T$'s are expressed as determinants involving
only the $T_{1,j,k}$'s. These particular Pl\"ucker relations are also known as the Desnanot-Jacobi
relation, used by Dodgson to devise his famous algorithm for the computation of determinants \cite{DOD}.
In \cite{RR}, this equation was slightly deformed by introducing a parameter $\lambda$ before
the second term on the r.h.s. and used to define the ``lambda-determinant", with a remarkable
expansion on alternating sign matrices, generalizing the usual determinant expansion over permutations.
Here we will not consider such a deformation, although we believe our constructions can be adapted to include this case as well (see \cite{SPY} for a general discussion, which however does not cover the $A_r$ case). 

Viewing the system \eqref{Tsys} as a three-term recursion relation in $k\in \Z$, 
it is clear that the solution is entirely
determined in terms of some initial data that covers two consecutive values of $k$, say $k=0,1$
and say all $j\in \Z$. In \cite{DFK09a},
an explicit expression for $T_{\al,j,k}$ was derived as a function of the initial data 
$\bx_0=\{T_{\al,j,0},T_{\al,j,1}\}_{\al\in I_r,j\in\Z}$. 
It involved expressing first $T_{1,j,k}$ as the partition function
for weighted paths on some particular target graph, with weights that are monomials 
of the initial data, and then
interpreting  $T_{\al,j,k}$ as the partition of non-intersecting families of such paths.
This interpretation was then extended to other initial data of the form
\begin{equation}\label{homodata}
\bx_{\bf k}=\{T_{\al,j,k_\al},T_{\al,j,k_\al+1}\}_{\al\in I_r,j\in\Z}
\end{equation}
where 
${\bf k}=(k_1,k_2,...,k_r)\in\Z^r$ is a Motzkin path of
length $r-1$, namely $k_{\al+1}-k_\al\in \{0,1,-1\}$ for all $\al=1,2,...,r-1$. 
In this construction, for each Motzkin path $\bf k$, the expressions for the $T_{\al,j,k}$
in terms of the initial data $\bx_{\bf k}$ are also partition functions of weighted 
paths on some target graph $\Gamma_{\bf k}$.

The equation \eqref{Tsys} is also connected to cluster algebras. 
In Ref. \cite{DFK08}, it was shown that the initial data sets 
$\bx_{\bf k}$ form a particular subset of clusters in a
suitably defined cluster algebra.
Roughly speaking, a cluster algebra \cite{FZI} is a dynamical system expressing the evolution of some
initial data set (cluster), with the built-in property that any evolved data is expressible as Laurent polynomials
of any other data set. This Laurent property or Laurent phenomenon turned out to be even more powerful
than expected, as all the known examples show that these polynomials have 
{\it non-negative integer} coefficients. The positivity conjecture of \cite{FZI} states that this
property holds in general. As an example, the above-mentioned lambda-determinant relation may be
viewed as an evolution equation in the same cluster algebra as in \cite{DFK08},
with $\lambda$ as a coefficient: the 
existence of an expansion formula of the lambda-determinant on alternating sign matrices is a 
manifestation of the positive Laurent phenomenon.
As another example,  the explicit expressions of \cite{DFK09a} for the
solutions $T_{\al,j,k}$ of the $A_r$ $T$-system as partition functions 
for positively weighted paths gives a direct proof of Laurent positivity 
for the relevant clusters. 

However, the set of initial data $\bx_{\bf k}$ \eqref{homodata} covered in \cite{DFK09a} is limited to sets of 
$T_{\al,j,k}$'s with fixed values of $k=k_\al$ independently of $j$. 
The most general set of initial data
should also allow for inhomogeneities in $j$, namely values of $k=k_{\al,j}$ varying with $j$ as well.
It is easy to see that the most general boundary condition consists in assigning fixed positive values 
$(a_{\al,j})_{\al\in I_r;j\in\Z}$ to $T_{\al,j,k_{\al,j}}$ along a ``stepped surface" (also called solid-on-solid
interface in the physics literature), namely such that $|k_{\al+1,j}-k_{\al,j}|=1$ and $|k_{\al,j+1}-k_{\al,j}|=1$
for all $\al\in I_r$ and $j\in \Z$.

In this paper we address the most general case of initial data for the $A_r$
$T$-system (\ref{Tsys}-\ref{bounT}). As we will show, initial data are in bijection with configurations
of the six-vertex model with face labels on a strip of square lattice of height $r-1$ and infinite width. 
For any such given set of initial data, we
derive an explicit expression for the solution $T_{\al,j,k}$ (\ref{Tsys}-\ref{bounT}) as the
partition function for $\al$ non-intersecting paths on a suitable network, in the spirit of Refs.
\cite{FZposit} and \cite{LP}, and with step weights that are Laurent monomials of the initial data. 
This completes the proof of the Laurent 
positivity of the solutions of the $A_r$ $T$-system for arbitrary initial data.

The paper is organized as follows. 

Our construction was originally inspired by Ref.\cite{FRISES}
which basically deals with the case of $A_1$ under the name of ``frises": 
the latter is reviewed in Section 2, where we make in particular
the connection between the ``frise" language and the solutions of the $A_1$ $T$-system
with arbitrary boundary data. Roughly speaking, the solution is expressed as the 
element of a matrix product taken along the boundary.

A warmup generalization to the case of $A_2$ is presented in Section 3, 
with the main Theorem \ref{genA2} giving an explicit solution for arbitrary boundary data, also as
an element of a matrix product taken along the boundary.

Section 4 is devoted to the general $A_r$ case. Starting from the path solution of \cite{DFK09a}
for some particular initial data, we construct various transfer matrices associated to the boundary,
with simple transformations under local elementary changes of the boundary (mutations).
For convenience, boundaries are expressed as configurations of the six-vertex model in an
infinite strip of finite height $r-1$.
These in turn encode a network, entirely determined by the boundary data.
The final result is an explicit formula Theorem \ref{maintheo} for the solution of the
$A_r$ $T$-system as the partition function for families of non-intersecting paths.

In Section 5, we study the restrictions of our results to the $Q$-system.

A few concluding remarks are gathered in Section 6.

\section{$A_1$ $T$-system and Frises}

In this section, we first review the results of \cite{FRISES}, and then rephrase them in terms of 
solutions to the $A_1$ $T$-system for arbitrary boundary conditions.

\subsection{Frises}

\subsubsection{Frise equation}
The frise equation reads\footnote{Our convention corresponds to $b\to -b$
in those of Ref \cite{FRISES}.}:
\begin{equation}\label{frise}
u_{a+1,b-1}u_{a,b}=1+u_{a+1,b}u_{a,b-1}
\end{equation}
for $a,b\in\Z$.

\subsubsection{Boundaries}

The most general (infinite) boundary condition is along 
a ``staircase", made of horizontal (h) and vertical (v) steps
of the form $h: (x,y)\to (x+1,y)$ and $v:(x,y)\to (x,y+1)$, giving rise to a
sequence of vertices $(x_j,y_j)$, $j\in \Z$. To each
vertex of the sequence we attach a positive number $a_j$, $j\in \Z$, and the boundary
condition for the system \eqref{frise} reads:
\begin{equation}\label{boundaone}
u_{x_j,y_j}=a_j \qquad (j\in \Z)
\end{equation}

The simplest such boundary is the sequence $...hvhvhv...$,
say with variable $a_{2x}$ at vertex $(x,x)$ and $a_{2x+1}$ at
vertex $(x,x+1)$, $x\in \Z$. We refer to it as the {\it basic staircase
boundary}.

The problem is now to find the solution $u_{a,b}$ of \eqref{frise}
with the boundary condition \eqref{boundaone}.


\subsubsection{Projection of $(x,y)$ on the boundary and step matrices}

\begin{figure}
\centering
\includegraphics[width=8.cm]{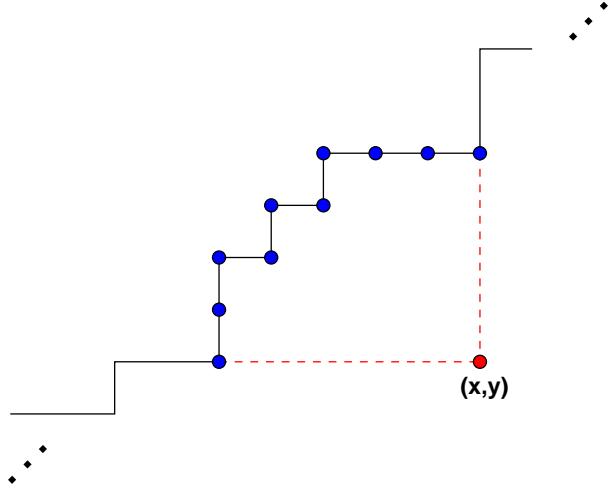}
\caption{\small A typical boundary for the frise and the projection of a point $(x,y)$
onto it. The vertices of the projection on the boundary are represented with blue circles.
Here, the word $w(x,y)$ reads $v^2hvhvh^3$.}
\label{fig:boundA1}
\end{figure}

The general solution at a point $(x,y)$ to the right of the boundary is expressed
solely in terms of the values \eqref{boundaone} taken by $u$
along the ``projection" of $(x,y)$ onto the boundary, defined as follows.

\begin{defn}[Projection]
\label{projaone}
The projection of $(x,y)$ onto the boundary $\{(x_j,y_j)\}_{j\in \Z}$ is the
sequence $(x_j,y_j)$, $j=t,t+1,...,t'$, where
$y_t=y$, $x_{t'}=x$, and the first step $t\to t+1$ is vertical,
while the last step $t'-1\to t'$ is horizontal.
\end{defn}

This is illustrated in Fig.\ref{fig:boundA1}. 
Alternatively the projection of $(x,y)$  is coded by the word 
$w(x,y)=$v...h
of length $t'-t$ with letters h and v, starting with v and ending with h
and coding the succesion of horizontal (h) and vertical (v) steps along the 
boundary between $(x_t,y_t)$ and $(x_{t'},y_{t'})$.
We also define the corresponding
sequence of boundary weights $\ba(x,y)=(a_t,a_{t+1},...,a_{t'})$.

\begin{defn}[Step matrices]\label{HVmat}
We define the two horizontal and vertical matrices
\begin{equation}
H(a,b)={1\over b} \begin{pmatrix} b & 0 \\ 1 & a \end{pmatrix}
\qquad  V(a,b)={1\over b}
\begin{pmatrix} a & 1 \\ 0 & b\end{pmatrix}
\end{equation}
\end{defn}

\subsubsection{Solution}

Given some boundary conditions, we associate
to the word $w$ and the
sequence $\ba$ the following $2\times 2$ matrix
product
\begin{equation}\label{mmatrix}
M(w,\ba)=V(a_t,a_{t+1})\cdots H(a_{t'-1},a_{t'})
\end{equation}
where the product extends over all the intermediate steps $i\to i+1$ between 
$t$ and $t'$ as coded by $w$, and involves the matrix $H(a_i,a_{i+1})$ if the step $i\to i+1$ is h
and $V(a_i,a_{i+1})$ if it is v. The result of \cite{FRISES} takes the following form:

\begin{thm}[\cite{FRISES}]\label{frisetheo}
The solution of \eqref{frise} subject to the boundary condition \eqref{boundaone} reads:
\begin{equation}\label{frisol}
u_{x,y}=a_{t'} \left(M(w(x,y),\ba(x,y))\right)_{1,1}
\end{equation}
with $M$ as in \eqref{mmatrix}.
\end{thm}

All matrices $V,H$ having elements that are positive Laurent monomials of the initial data,
the general Laurent positivity of the solution follows:

\begin{cor}
The general solution of \eqref{frise} subject to the boundary condition \eqref{boundaone}
is a Laurent polynomial of its initial data $\{a_j\}_{j\in \Z}$, with non-negative integer coefficients.
\end{cor}

\begin{example}[The basic staircase boundary]

For any $x>y\in \Z_{\geq 0}$, we have a projection on the boundary
with $w(x,y)=(vh)^{x-y}$, and $\ba(x,y)=(a_{2y},a_{2y+1},...,a_{2x})$.
We deduce that
\begin{equation}\label{matstair}
M(w(x,y),\ba(x,y))=\prod_{i=y}^{x-1} V(a_{2i},a_{2i+1})H(a_{2i+1},a_{2i+2})
\end{equation}
and the solution reads:
\begin{equation}\label{restair}
u_{x,y}=a_{2x}\, \left(M(w(x,y),\ba(x,y))\right)_{1,1}
\end{equation}
Explicitly, we compute the two-step matrix:
\begin{equation}\label{vh} 
M(a,b,c)=V(a,b)H(b,c)= 
{1\over c} \begin{pmatrix} {a c+1\over b} & 1\\ 1 & b \end{pmatrix} 
\end{equation}

\end{example}

\subsubsection{Mutations}

Note that we may move
from one boundary to another by elementary ``mutations"\footnote{The 
term ``mutation" is borrowed from cluster algebras, as this elementary move indeed corresponds
to a mutation in the associated cluster algebra of \cite{DFK09a}.}, namely the 
local substitution $(v,h)\rightarrow (h,v)$ on the boundary (forward mutation)
or $(h,v)\rightarrow (v,h)$ (backward mutation), while the 
sequence $\ba$ is updated using the frise relation \eqref{frise}. In particular, we may in principle
reach any boundary from the basic staircase one, by possibly infinitely many 
such mutations. 

The effect of such a mutation is easily obtained by computing
the corresponding matrix transformation within $M(w,\ba)$.
It basically corresponds to the following identity:

\begin{lemma}\label{mutaA1}
For all $a,b,c>0$, we have:
\begin{equation}\label{mutamat}
\quad\quad \raisebox{-.7cm}{\hbox{\epsfxsize=2.5cm \epsfbox{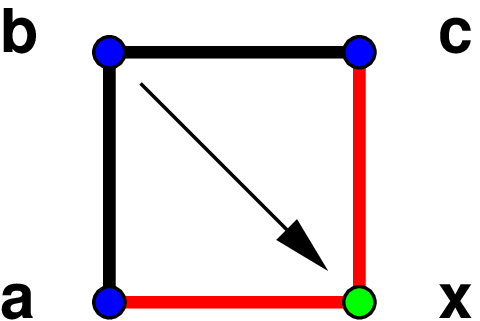}}}\quad
V(a,b) H(b,c) =H(a,x)V(x,c), \qquad x={1+a c\over b}
\end{equation}
\end{lemma}

This may be understood as a matrix representation of the mutation via 
the commutation of the matrices $V$ and $H$, which acquire the new 
boundary value $x$ in replacement for $b$. 
This mutation affects all values of $u_{m,p}$ such that the projection of $(m,p)$
contains the new boundary point with value $x$. 

We may deduce the general formula \eqref{frisol} from that for the basic staircase
boundary, by induction under mutation. In general a mutation simply switches two
consecutive matrices $VH\to HV$ in the product $M$. We must be careful with 
mutations that update the extremal vertices of the projection of $(x,y)$, namely
in the two cases:
(i) when  $w$ starts with vh, updated into hv or (ii) when $w$ ends up with vh,
updated into hv. 
We note however that
$H(a,b)_{1,j}=bV(a,b)_{j,1}=\delta_{i,j}$ hence in the updated matrix $M$
we may: (i) drop the first matrix factor $H$ (ii) drop the last matrix factor $V$,
but replace the scalar prefactor by the new updated vertex value, and the formula 
\eqref{frisol} follows.

\subsection{The $A_1$ T-system}

\subsubsection{T-system}

The $A_1$ T-system reads:
\begin{equation}\label{tsysa1}
T_{j,k+1}T_{j,k-1}=T_{j+1,k}T_{j-1,k}+1
\end{equation}
where we use the shorthand notation $T_{j,k}=T_{1,j,k}$
for $j,k\in \Z$. Note that this splits into two independent systems for fixed value
of $j+k$ modulo 2.

In the case when $j+k=0$ modulo 2,
we immediately see that changing to ``light cone" coordinates: 
$a={j+k\over 2}$ and $b={j-k\over 2}$, we have that $u_{a,b}=T_{j,k-1}$
satisfies the frise equation \eqref{frise}. So the two problems are equivalent.
Analogously, when $j+k=1$ modulo 2, we take $a={j+k-1\over 2}$ and $b={j-k+1\over 2}$
and $u_{a,b}=T_{j,k}$.

\subsubsection{Boundaries}

The fundamental boundary for the T-system is obtained by fixing the
values of say $T_{2j+1,0}$ and $T_{2j,1}$ for all $j\in \Z$. It corresponds to the
basic staircase boundary in the case $j+k=1$ modulo 2 above, with 
$T_{2j,1}=a_{2j}$ and $T_{2j+1,0}=a_{2j+1}$.

Other boundaries are mapped in an obvious manner. 

\subsubsection{Path solution}

In \cite{DFK09a}, an explicit path formulation was derived for the solution $T_{j,k}$ 
for the fundamental boundary condition. Defining
the $4\times 4$ transfer matrix
$$ \cT(u,v,w)= \begin{pmatrix} 0 & 1 & 0 & 0\\ u & 0 & 1 & 0\\ 0 & v & 0 & 1\\
0 & 0 & w & 0 \end{pmatrix} $$
we have

\begin{thm}[\cite{DFK09a}]
The solution of the $A_1$ $T$-system \eqref{tsysa1} for the fundamental 
boundary condition with initial data $\{ T_{2j+1,0},T_{2j,1}\}_{j\in\Z}$ reads for $j+k=1$ mod $2$:
\begin{eqnarray}
T_{j,k}&=& T_{j+k,0} \left(\prod_{i=j-k}^{j+k-1}\cT(u_i,v_i,w_i) \right)_{1,1}\nonumber \\
u_i&=& {T_{i,1}\over T_{i+1,0}}, \quad v_i={1\over T_{i,0}T_{i+1,1}},
\quad w_i={1\over u_i} \label{tsysol}
\end{eqnarray}
\end{thm}

Here the matrix $\cT(i,i+1)\equiv \cT(u_i,v_i,w_i)$ is interpreted as the transfer matrix 
from time $i$ to time $i+1$
for weighted paths with steps $a\to a\pm 1$ on the integer segment $[0,3]$, 
with time-dependent step weights
$0\to 1: u_i$, $1\to 2: v_i$ and $2\to 3: w_i$, the other weights being equal to $1$.

\subsubsection{Gauge invariance}
The above formula \eqref{tsysol} remains clearly unchanged if we transform the matrix
$\cT$ into the following: $\cT(i,i+1) \to {\tilde \cT}(i,i+1)=L_i \cT(i,i+1) L_{i+1}^{-1}$,
for any invertible matrix $L_i$ such that $(L_i)_{1,j}=\delta_{j,1}$
and $(L_i)_{j,1}=\delta_{j,1}$.

To make the contact with the frise solution, let us define ${\tilde \cT}$ as above,
by use of the matrix
$$ L_i= \begin{pmatrix} 1 & 0 & 0 & 0 \\ 0 & 1 & 0 & 0 \\ 0 & 0 & T_{i,0} & 0 \\
0 & 0 & 0 & T_{i, 1}\end{pmatrix} $$

\subsubsection{Comparison with the frise solution}

Let us now compute the ``two-step" transfer matrix at times $i,i+1$ for $i=j+k$ modulo 2:
$ {\tilde \cT}(i,i+2)= {\tilde \cT}(i,i+1)  {\tilde \cT}(i+1,i+2)$,
with the weights as in \eqref{tsysol}, namely:
\begin{eqnarray*}
&&u_{i}={b_{i}\over b_{i+1}},\quad 
v_{i}={1\over a_{i}a_{i+1}},\quad w_{i}={1\over u_{i}} \\
&&u_{i+1}={a_{i+1}\over a_{i+2}},\quad 
v_{i+1}={1\over b_{i+1}b_{i+2}},\quad w_{i+1}={1\over u_{i+1}}
\end{eqnarray*}
where we have introduced $a_i=T_{i,0}$, $a_{i+1}=T_{i+1,1}$, $a_{i+2}=T_{i+2,0}$,
$b_i=T_{i,1}$, $b_{i+1}=T_{i+1,0}$,  and $b_{i+2}=T_{i+2,1}$, while
$L_i={\rm diag}(1,1,a_i,b_i)$, $L_{i+1}={\rm diag}(1,1,b_{i+1},a_{i+1})$
and $L_{i+2}={\rm diag}(1,1,a_{i+2},b_{i+2})$.

We get:
$$  {\tilde \cT}(i,i+2)=
\begin{pmatrix}
{a_{i+1}\over a_{i+2}} & 0 & {1\over a_{i+2}} & 0\\
0 & {1+b_{i}b_{i+2}\over b_{i+1} b_{i+2}} & 0 & {1\over b_{i+2}}\\
{1\over a_{i+2}}  & 0 & {1+a_{i}a_{i+2}\over a_{i+1}a_{i+2}} & 0\\
0 & {1\over b_{i+2}} & 0 & {b_{i+1}\over b_{i+2}}
\end{pmatrix} 
$$

This $4\times 4$ matrix clearly decomposes into two independent $2\times 2$
linear operators acting respectively on components $1,3$ and $2,4$.
The corresponding matrices are respectively:
\begin{eqnarray*}
{\tilde \cP}(i,i+2)&=& \begin{pmatrix}
{a_{i+1}\over a_{i+2}} & {1\over a_{i+2}} \\
{1\over a_{i+2}} & {1+a_{i}a_{i+2}\over a_{i+1}a_{i+2}}
\end{pmatrix}=H(a_{i},a_{i+1})V(a_{i+1},a_{i+2})\\
{\tilde \cQ}(i,i+2)&=&\begin{pmatrix}
{1+b_{i}b_{i+2}\over b_{i+1} b_{i+2}} & {1\over b_{i+2}}\\
{1\over b_{i+2}} &  {b_{i+1}\over b_{i+2}}
\end{pmatrix}=V(b_i,b_{i+1})H(b_{i+1},b_{i+2})
\end{eqnarray*}

We may now use the gauge-transformed and
reduced two-step transfer matrix ${\tilde \cP}(i,i+2)$
instead of ${\cT}$ in \eqref{tsysol}.
Indeed, in the product over steps from $j-k$ to $j+k-1$,
we may pair up consecutive ${\cT}$ matrices in \eqref{tsysol} to express 
it in terms of the $\cP$'s, and then substitute the latter with the $\tilde \cP$'s, leading to:
$$T_{j,k}=T_{j+k,0} \left( 
\prod_{i=0}^{k-1} {\tilde \cP}(j-k+2i,j-k+2i+2) \right)_{1,1}$$
Noting moreover that $H(a,b)_{1,j}=b V(a,b)_{j,1}=\delta_{j,1}$, we may
rewrite this as
\begin{equation}\label{upsidown}
T_{j,k}=T_{j+k-1,1} \left( 
\prod_{i=0}^{k-2} V(a_{j-k+2i+1},a_{j-k+2i+2})H(a_{j-k+2i+2},a_{j-k+2i+3}) \right)_{1,1}
\end{equation}
Assuming that $j+k=1$ modulo 2, we see that in light-cone coordinates with
$x={j+k-1\over 2}$ and $y={j-k+1\over 2}$, equation \eqref{upsidown}
amounts to equations (\ref{matstair}-\ref{restair}), as the projection
of $(x,y)$ on the boundary staircase starts at $t=j-k+1$ and ends at $t'=j+k-1$.

\subsubsection{Mutations and arbitrary boundary}

As in the frise case, this identification gives us access to mutations, 
via the $VH\leftrightarrow HV$ identity \eqref{mutamat}. Starting
from \eqref{upsidown}, we may iteratively apply forward/backward mutations
to the basic staircase boundary to get any other boundary (up to global
translations) of the form $\{ T_{j,k_j}\}_{j\in \Z}$ with a sequence $k_j\in\Z$
such that $|k_{j+1}-k_j|=1$. 
Let us denote by $(j_0,k_{j_0})$ and $(j_1,k_{j_1})$ the extremities
of the projection of $(j,k)$ onto the boundary, namely such that $j_0-k_{j_0}=j-k$,
$j_1+k_{j_1}=j+k$, $j_0$ maximal and $j_1$ minimal.

We deduce that the general solution for arbitrary staircase boundary reads:
$$T_{j,k}=T_{j_1,k_{j_1}}
\Big(V(T_{j_0,k_{j_0}},T_{j_0+1,k_{j_0+1}})...
H(T_{j_1-1,k_{j_1-1}},T_{j_1,k_{j_1}})\Big)_{1,1}
$$
where the product is taken along the projection of $(j,k)$ on the boundary,
with a matrix $V$ per vertical step and $H$ per horizontal step.

\section{The $A_2$ T-system with arbitrary boundary}

Before going to the general $A_r$ case, we derive the $A_2$ solution in detail.

\subsection{T-system}

The $A_2$ T-system reads:
\begin{eqnarray}
T_{1,j,k+1}T_{1,j,k-1}&=&T_{1,j+1,k}T_{1,j-1,k}+T_{2,j,k}\nonumber \\
T_{2,j,k+1}T_{2,j,k-1}&=&T_{2,j+1,k}T_{2,j-1,k}+T_{1,j,k}
\label{tsysa2}
\end{eqnarray}
for $j,k\in \Z$. Note that this splits again into two independent systems for $T_{\al,j,k}$
with fixed value of $\al+j+k$ modulo 2. These indices run over two consecutive layers
of the centered cubic lattice $\al=1$ and $\al=2$, which form two square lattices, 
the vertices of the second layer lying at the vertical of the centers of the faces of the first layer.

\subsection{Boundaries}

\begin{figure}
\centering
\includegraphics[width=12.cm]{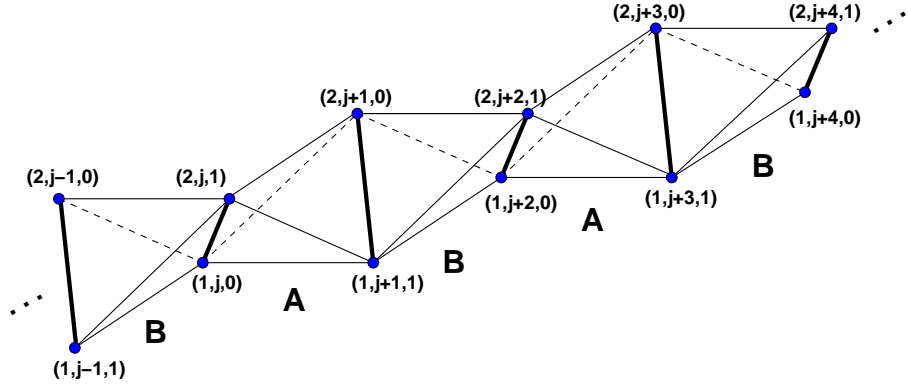}
\caption{\small The basic staircase boundary for the $A_2$ $T$-system. 
We have indicated the corresponding succession of edges (thick black line)
and the two types of tetrahedrons A,B that connect them.}
\label{fig:stepa2}
\end{figure}

The fundamental boundary considered in Ref. \cite{DFK09a} involves fixing the values
of the $T_{\al,j,k}$
with $\al=1,2$ $k=0,1$, and $j\in \Z$, with fixed parity of $\al+j+k$ (say even). 
We refer to this boundary as the {\it basic staircase} boundary, in reference to the 
$A_1$ case.
It can be viewed as an infinite strip made of a succession of four kinds of vertices
(see Fig.\ref{fig:stepa2}). 
We may also view this strip as a succession of edges of the form
$e_j=(1,j,0)-(2,j,1)$, $f_{j+1}=(1,j+1,1)-(2,j+1,0)$, 
$e_{j+2}=(1,j+2,0)-(2,j+2,1)$, etc. for $j-1\in 2\Z$ 
(thick black lines in Fig.\ref{fig:stepa2}). Two such
consecutive edges define a tetrahedron. 
The basic staircase may therefore
be viewed as the alternating succession of two kinds of tetrahedrons denoted by 
A (defined by $e_j,f_{j+1}$) and B (defined by $f_{j-1},e_j$).

\subsection{Solution for the basic staircase boundary}

In Ref. \cite{DFK09a}, the solution $T_{1,j,k}$ was expressed in terms of paths on a target
graph with 6 vertices and with time-dependent edge weights involving only the boundary values.
These weights are coded by a $6\times 6$ transfer matrix. Defining:
\begin{equation}
\cT(s,t,u,v,w)=\begin{pmatrix}
0 & 1 & 0  & 0 & 0 & 0 \\
s & 0 & 1 & 0 & 0 & 0 \\
0 & t & 0 & 1 & 1 & 0\\
0 & 0 & u & 0 & 0 & 0\\
0 & 0 & v & 0 & 0 & 1\\
0 & 0 & 0 & 0 & w & 0
\end{pmatrix}
\end{equation}
and using the notation
\begin{equation}
s_i={T_{1,i,1}\over T_{1,i+1,0}}, \ \  t_i={T_{2,i,1}\over T_{1,i,0}T_{1,i+1,1}}, \ \ 
u_i= {T_{1,i+1,0}T_{2,i-1,1}\over T_{1,i,1}T_{2,i,0}},\ \  
v_i={T_{1,i+1,0}\over T_{2,i,0}T_{2,i+1,1}},\ \ w_i={T_{2,i+1,0}\over T_{2,i,1}}
\end{equation}
we have:

\begin{thm}[\cite{DFK09a}]
The solution of the $A_2$ $T$-system for $\al=1$ reads:
\begin{equation}\label{forma2}
T_{1,j,k}=T_{1,j+k,0} \left( \prod_{i=j-k}^{j+k-1} \cT(s_i,t_i,u_i,v_i,w_i) \right)_{1,1}
\end{equation}
\end{thm}

\subsection{Reduced transfer matrix}

As before we note that the two-step transfer matrix 
$\cT(i,i+2)=\cT(i,i+1)\cT(i+1,i+2)$, with 
$\cT(i,i+1)\equiv \cT(s_i,t_i,u_i,v_i,w_i)$,
is again decomposable into two linear operators acting respectively
on components $(1,3,6)$ and $(2,4,5)$. Explicitly:
$$\cT(s,t,u,v,w)\cT(s',t',u',v',w')=\begin{pmatrix} 
s' & 0 & 1 & 0 & 0 & 0 \\
0 & s+t' & 0 & 1 & 1 & 0\\
t s' & 0 & t+u'+v' & 0 & 0 & 1\\
0 & ut' & 0 & u & u & 0 \\
0 & vt' & 0 & v & v+w' & 0 \\
0 & 0 & wv' & 0 & 0 & w
\end{pmatrix}$$
Defining:
$$ \cP(i,i+2)=\begin{pmatrix}
s_{i+1} & 1 & 0\\
t_i s_{i+1} & t_i+u_{i+1}+v_{i+1} & 1\\
0 & w_i v_{i+1}  & w_i
\end{pmatrix} $$
we may rewrite \eqref{forma2} as:
\begin{equation}\label{sola2}
T_{1,j,k}=T_{1,j+k,0} \Big( \prod_{i=0}^{k-1} 
\cP(j-k+2i,j-k+2i+2)\Big)_{1,1}
\end{equation}

\subsection{Gauge transformation, tetrahedron decomposition, and frise analogy}

As before, we note that any gauge transformation of the form
$\cP(i,i+2) \to {\tilde \cP}(i,i+2)=L_i \cP(i,i+2) L_{i+2}^{-1}$ and such that 
$(L_i)_{1,j}=(L_i)_{j,1}=\delta_{j,1}$ leaves the formula \eqref{sola2} 
invariant.

We choose 
$$L_i=\begin{pmatrix} 1 & 0 & 0\\
0 & T_{1,i,0} & 0\\
0 & 0 & T_{2,i,1} \end{pmatrix}$$
Defining 
$$M(a,b,c,u,v,w)=\begin{pmatrix} {b\over c} & {1\over c} & 0\\
{u\over c} & {u\over bc}+{au\over bv}+{a\over vw} & {a\over w}\\
0 & {1\over w} & {v\over w} \end{pmatrix}$$
we have 
$${\tilde \cP}(i,i+2)=
M(T_{1,i,0},T_{1,i+1,1},T_{1,i+2,0},T_{2,i,1},T_{2,i+1,0},T_{2,i+2,1})
$$

The matrix $M$ may be further decomposed as follows:
$$M(a,b,c,u,v,w)=A(a,b,u,v)B(b,c,v,w)$$
where
\begin{eqnarray*}
\quad \raisebox{-1.cm}{\hbox{\epsfxsize=2.5cm \epsfbox{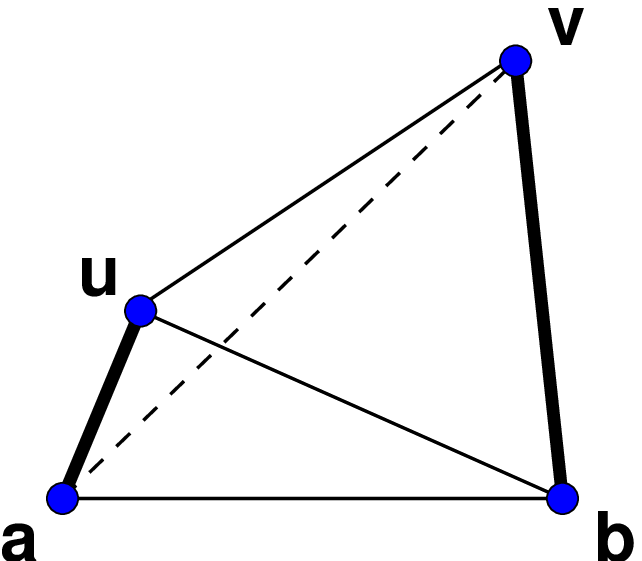}}}\qquad \qquad 
 A(a,b,u,v)&=&\begin{pmatrix}1 & 0 & 0 \\
{u\over b} & {a u \over b v} & {a\over v}\\
0 & 0 & 1 \end{pmatrix} \\ 
\quad \raisebox{-1.cm}{\hbox{\epsfxsize=2.5cm \epsfbox{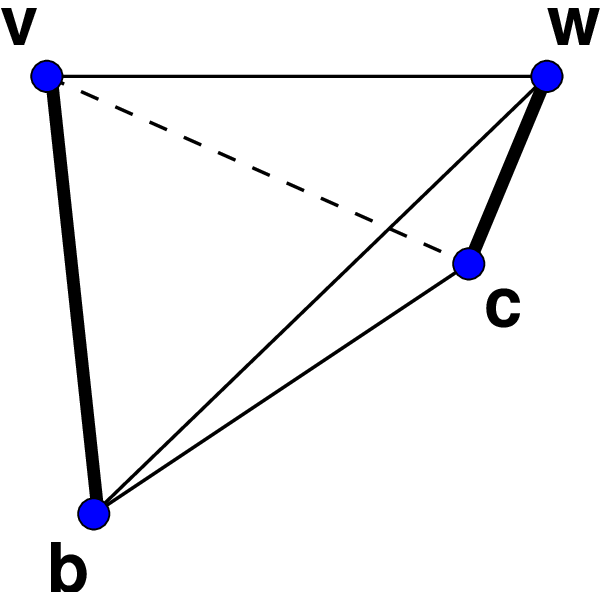}}}\qquad \qquad 
B(b,c,v,w)&=&\begin{pmatrix} {b\over c} & {1\over c} & 0\\
0 & 1 & 0\\
0 & {1\over w} & {v \over w} 
\end{pmatrix} 
\end{eqnarray*}

In view of our interpretation of the boundary (see Fig.\ref{fig:stepa2}),
the matrices $A$ and $B$ may be associated to the tetrahedrons A and B.
The arguments are the values of $T$ at the vertices of the tetrahedrons.
We write 
$$ A(i,i+1)=A(T_{1,i,1},T_{1,i+1,0},T_{2,i,0},T_{2,i+1,1})\quad
{\rm and} \quad B(i,i+1)=B(T_{1,i,0},T_{1,i+1,1},T_{2,i,1},T_{2,i+1,0})$$
and finally we may rewrite \eqref{sola2} as:
\begin{equation}\label{a2sol}
T_{1,j,k}=T_{1,j+k,0} \Big( \prod_{i=0}^{k-1} 
A(j-k+2i,j-k+2i+1)B(j-k+2i+1,j-k+2i+2)\Big)_{1,1}
\end{equation}

We may now interpret this result as a generalization of the frise result. 
By analogy with the frise solution,
let us define the projection of $(1,j,k)$ on the boundary as the portion of
the boundary between the edge $f_{j-k+1}$ and the edge $f_{j+k-1}$.
We have:

\begin{thm}
The general solution of the $A_2$ $T$-system with the basic staircase boundary
for $\al=1$ reads:
\begin{equation}\label{frisa2}
T_{1,j,k}=T_{1,j+k-1,1} \Big(\prod_{i=0}^{k-2} 
B(j-k+2i+1,j-k+2i+2) A(j-k+2i+2,j-k+2i+3)\Big)_{1,1}
\end{equation}
\end{thm}
\begin{proof}
We start from \eqref{a2sol} and use the fact that $A(a,b,u,v)_{1,j}=b B(a,b,u,v)_{j,1}=\delta_{j,1}$
to eliminate the first ($A$) and last ($B$) matrices in the product on the r.h.s.
\end{proof}

The product extends over the sequence of tetrahedrons along the
projection of $(1,j,k)$ onto the boundary. We may think of the two tetrahedron matrices 
$A$, $B$ as a generalization of the horizontal and vertical matrices of the $A_1$ case,
but more general boundaries involve four more such matrices, as discussed below.

\subsection{Other boundaries: two tetrahedrons and four parallelograms}

\begin{figure}
\centering
\includegraphics[width=9.5cm]{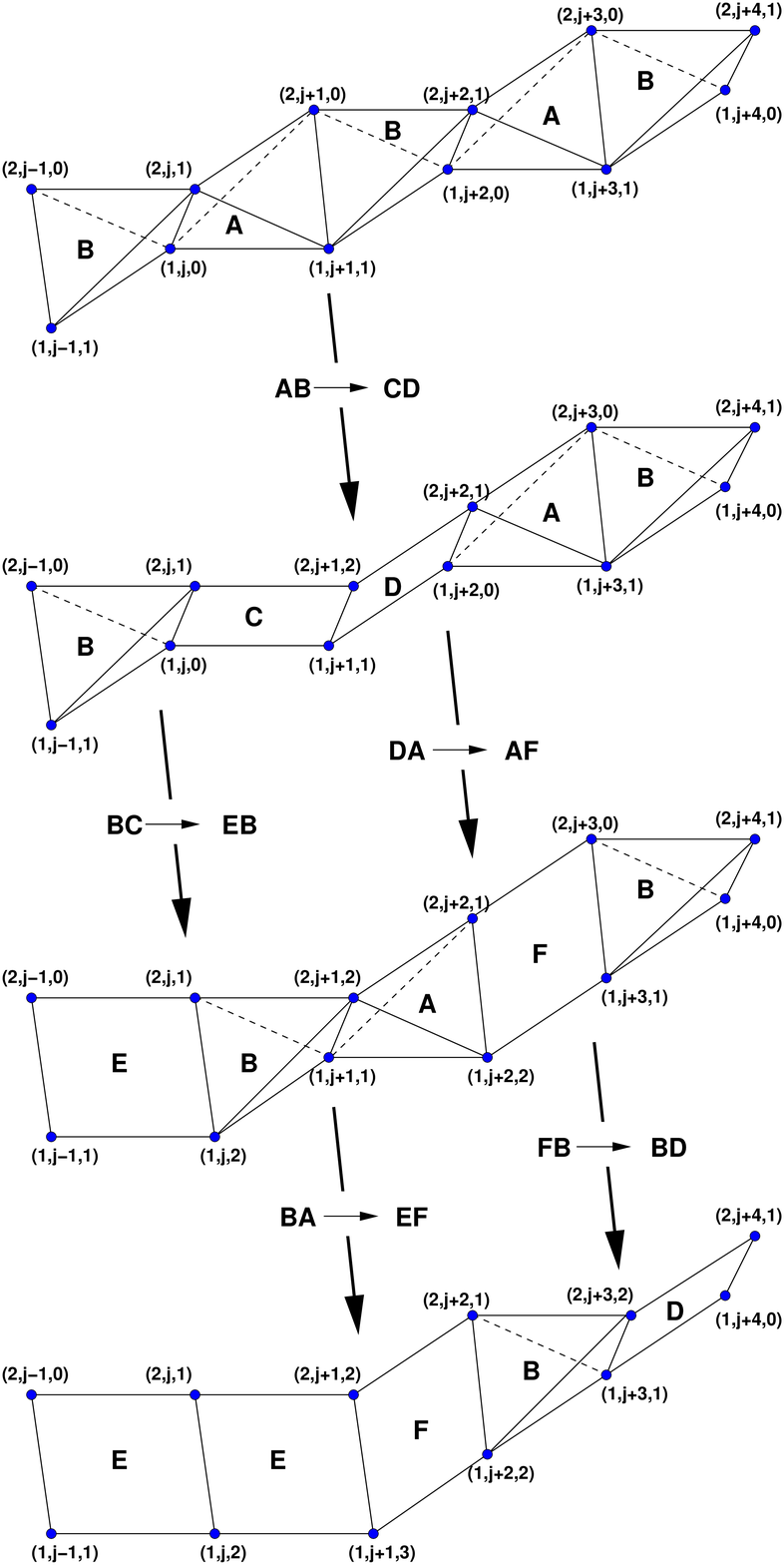}
\caption{\small Various mutations of the basic staircase boundary 
for the $A_2$ $T$-system. 
We have indicated the tetrahedrons A,B and the parallelograms C,D,E,F 
that may connect two consecutive edges of the boundary.}
\label{fig:muta2}
\end{figure}

The most general boundaries are obtained from the basic staircase
via local elementary moves (forward/backward  mutations) 
corresponding to one application
of one of the system relations. The effect is of flipping a single thick edge
of the boundary in the following manner: denoting by 
$e_{j,k}=(1,j,k)-(2,j,k+1)$ and $f_{j,k}=(1,j,k+1)-(2,j,k)$, we have
the two possible elementary (forward) moves:
\begin{eqnarray}
&&\ldots, e_{j,k},\ldots \to \ldots f_{j,k+1}\ldots 
\label{mut1}\\
&&\ldots, f_{j,k},\ldots \to \ldots e_{j,k+1}\ldots 
\label{mut2}
\end{eqnarray}
It is easy to see that this gives rise to six possible relative positions for two
consecutive edges:
$$ (e_{j,k},f_{j+1,k}), (e_{j,k},e_{j+1,k+1}), (e_{j,k},e_{j+1,k-1}),(f_{j,k},e_{j+1,k}),
(f_{j,k},f_{j+1,k+1}),(f_{j,k},f_{j+1,k-1}) $$
which define two tetrahedrons and four parallelograms, respectively denoted by
A,C,D,B,E,F (see Fig.\ref{fig:muta2}).

Note that the boundary is entirely specified by a Motzkin path and the edge at one of
its vertices. Indeed the transition from an edge to the next changes $k\to k, k+1$ or $k-1$,
and the nature of the edge ($e$ or $f$) is switched only if $k$ is unchanged. So keeping
a record of the variable $k$ is sufficient, and the record is a Motzkin path 
(or equivalently an infinite word in 3 letters).

\subsection{Triangle decomposition}

\begin{figure}
\centering
\includegraphics[width=12.cm]{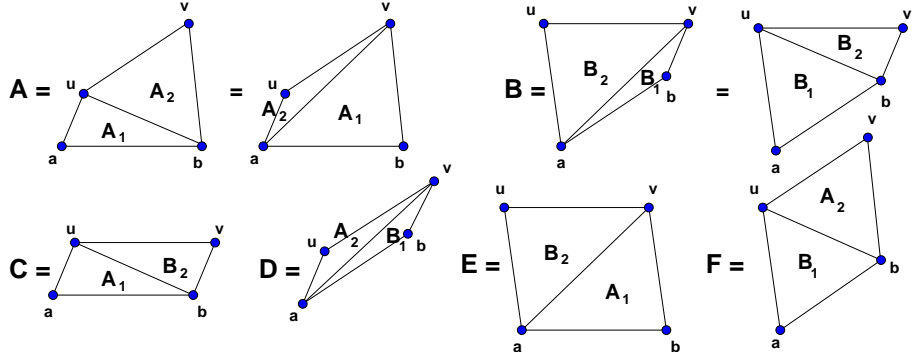}
\caption{\small Triangle decomposition A$_1$,A$_2$,B$_1$,B$_2$ of the 
tetrahedrons A,B and parallelograms 
C,D,E,F. We use the same letter for different triangles that share the same 
transfer matrices (see eq.\eqref{ABXY}).}
\label{fig:tetramata2}
\end{figure}

As already mentioned, we may further decompose the tetrahedrons A and
B, as well as the parallelograms C,D,E,F,  into pairs of triangles, as indicated
in Fig.\ref{fig:tetramata2}. 
To all triangles labelled A$_1$,A$_2$,B$_1$,B$_2$, we associate the following $3\times 3$ 
``triangle" matrices with 3 parameters equal to the values 
of $T_{\al,j,k}$ at their vertices.
We have:
\begin{eqnarray}
A_1(a,b,u)&=&\begin{pmatrix}
1 & 0 & 0 \\
{u\over b} & {a \over b} & 0 \\
0 & 0 & 1
\end{pmatrix}
\quad 
A_2(b,u,v)= \begin{pmatrix}
1 & 0 & 0 \\
0 & {u\over v} & {b \over v} \\
0 & 0 & 1
\end{pmatrix}\nonumber \\
B_1(a,b,u)&=&\begin{pmatrix}
{a\over b} & {1\over b} & 0 \\
0 & 1 & 0 \\
0 & 0 & 1
\end{pmatrix}
\quad 
B_2(b,u,v)=\begin{pmatrix}
1 & 0 & 0\\
0 & 1 & 0\\
0 & {1\over v} & {u\over v} \\
\end{pmatrix}\label{ABXY}
\end{eqnarray}
The two tetrahedrons A,B, correspond to the matrices
\begin{eqnarray}A(a,b,u,v)&=&A_1(a,b,u)A_2(b,u,v)=A_2(a,u,v)A_1(a,b,v)\nonumber \\
B(a,b,u,v)&=&B_1(a,b,u)B_2(b,u,v)=B_2(a,u,v)B_1(a,b,v)\label{ambigless}
\end{eqnarray}
independent of the two triangle decompositions. The parallelograms C,D,E,F
of Fig.\ref{fig:tetramata2} have unique triangle decompositions, to which we attach 
the following matrices:
\begin{eqnarray*}C(a,b,u,v)&=&A_1(a,b,u)B_2(b,u,v)={\small \begin{pmatrix}
1 & 0 & 0\\
{u\over b} & {a\over b} & 0\\
0 & {1\over v} & {u\over v}
\end{pmatrix}}\\
 D(a,b,u,v)&=&A_2(a,u,v)B_1(a,b,v)=\begin{pmatrix}
{a\over b} & {1\over b} & 0\\
0 & {u\over v} & {a\over v}\\
0 & 0 & 1
\end{pmatrix}\\
E(a,b,u,v)&=&B_2(a,u,v)A_1(a,b,v)=\begin{pmatrix}
1 & 0 & 0\\
{v\over b} & {a\over b} & 0\\
{1\over b} & {a\over v b} & {u\over v}
\end{pmatrix}\\ 
F(a,b,u,v)&=&B_1(a,b,u)A_2(b,u,v)=\begin{pmatrix}
{a\over b} & {u\over b v} & {1\over v}\\
0 & {u\over v} & {b\over v}\\
0 & 0 & 1
\end{pmatrix}
\end{eqnarray*}
Note that the products are taken in a specific order, 
namely that the matrix for the triangle which lies on the left 
of the diagonal of the parallelogram multiplies that on the 
right from the left.

\subsection{Mutations via triangles and the general formula}

\begin{figure}
\centering
\includegraphics[width=12.cm]{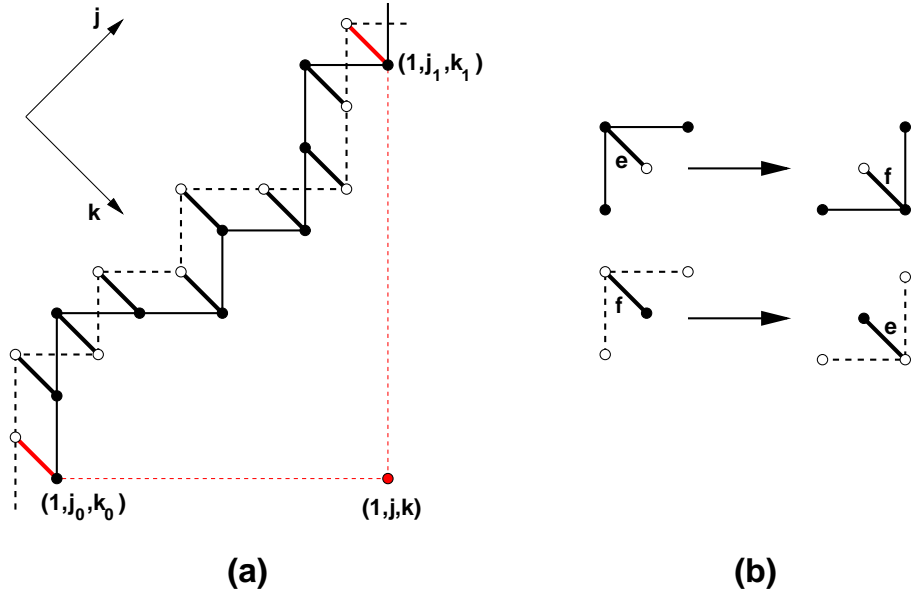}
\caption{\small The projection of a point $(1,j,k)$ onto a
typical boundary for the $A_2$ T-system (a), viewed in vertical projection
onto the bottom plane. The vertices of the bottom (rep. top) 
layers are represented as filled (resp. empty) circles. The boundary edges
are represented as thick diagonal lines.
We have indicated (b) the two types of forward mutations corresponding to $e\to f$ and $f\to e$
as in (\ref{mut1}-\ref{mut2}).}
\label{fig:bounda2}
\end{figure}

The general formula for $T_{1,j,k}$ for an arbitrary boundary reads as follows.
First, we may view the boundary above in yet another manner, by projecting 
it vertically onto the bottom plane, as illustrated in Fig. \ref{fig:bounda2}. It is the
superimposition of the two sets of boundary vertices in the bottom (resp. top) layers
(represented as filled (resp. empty) circles),
which are both staircase boundaries of the type considered in the $A_1$ case,
drawn on the two corresponding shifted square lattice layers in thick (resp. dashed)
lines. These two staircases are constrained by the condition that their vertices must be connected
via $e$ or $f$ edges only (diagonal thick black lines in Fig.\ref{fig:bounda2}).

\begin{defn}
The projection of the point $(1,j,k)$ onto the boundary is the portion of 
boundary (finite sequence of edges) between the 
edges containing $(1,j_0,k_0)$ and $(1,j_1,k_1)$, repectively such that 
$j_0-k_0=j-k$ and $j_1+k_1=j+k$ with $j_0$ maximal and $j_1$ minimal.
\end{defn}

This coincides with the definition for the $A_1$ T-system, using the staircase 
of the bottom layer only.
The corresponding finite sequence of edges corresponds alternatively
to a finite sequence of triangles according to the decomposition above,
modulo the two-fold ambiguities of decomposition of the tetrahedrons A and B. 

To this sequence
we associate the matrix:
\begin{equation}\label{Mmata2}
M(j,k)=\prod_{{\rm triangle}\ {\rm Z}\atop
{\rm w/vertex \ values}\ x,y,z } Z(x,y,z) 
\end{equation}
where for each triangle 
Z=A$_1$,A$_2$,B$_1$,B$_2$,
along the sequence
we multiply by the corresponding triangle matrix $Z=A_1,A_2,B_1,B_2$.
Note that, due to the identities \eqref{ambigless}, this definition is independent
of the particular choice of triangle decomposition of the possible tetrahedrons
along the boundary.
We have:

\begin{thm}\label{genA2}
The solution of the $A_2$ $T$-system for arbitrary boundary reads for $\al=1$:
\begin{equation}\label{gena2}
T_{1,j,k}=T_{1,j_1,k_1}\, M(j,k)_{1,1} 
\end{equation}
with $M(j,k)$ as in \eqref{Mmata2}, with the product extending 
over the projection of $(1,j,k)$ onto the boundary.
\end{thm}

Before proving the Theorem by induction under mutation, 
let us describe the mutations of the boundary in more detail.
The two possible mutations (\ref{mut1}-\ref{mut2}) correspond
to a local transformation of the chain of triangles that forms the boundary,
namely it replaces a pair of adjacent triangles sharing the initial boundary edge
with a new pair of adjacent triangles sharing the mutated boundary edge.
Using the definition \eqref{ABXY}, we get the following Lemma, 
generalizing Lemma \ref{mutaA1}:

\begin{lemma}\label{boundamut}
For all $a,b,c,u,v,w>0$ we have:
\begin{eqnarray}
\quad\quad \raisebox{-1.cm}{\hbox{\epsfxsize=3.5cm \epsfbox{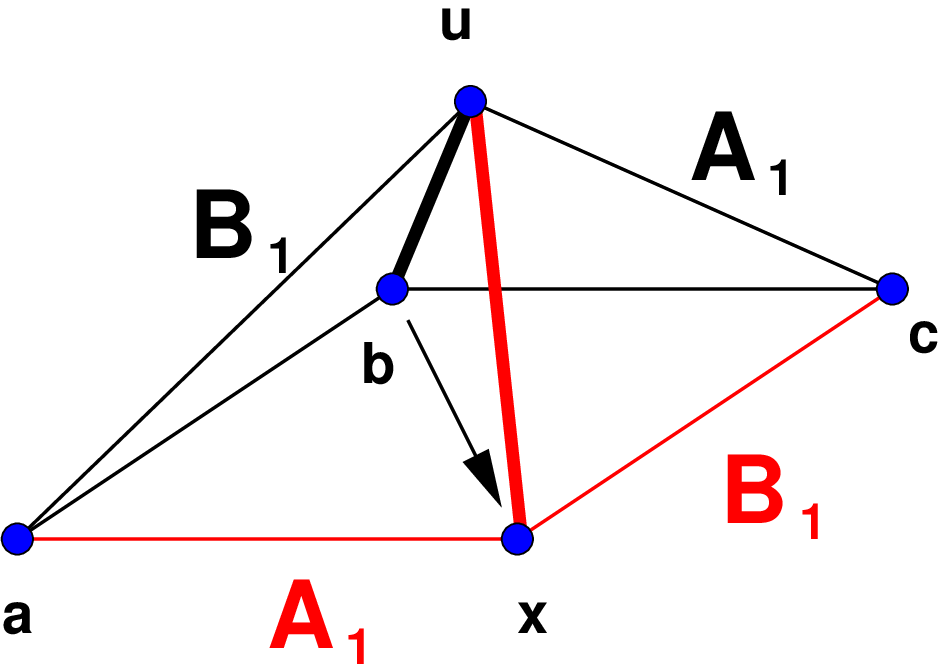}}}
\quad B_1(a,b,u)A_1(b,c,u)&=&A_1(a,x,u)B_1(x,c,u), 
\quad x={a c+u\over b}\label{mutone}\\
\quad\quad \raisebox{-1.cm}{\hbox{\epsfxsize=3.5cm \epsfbox{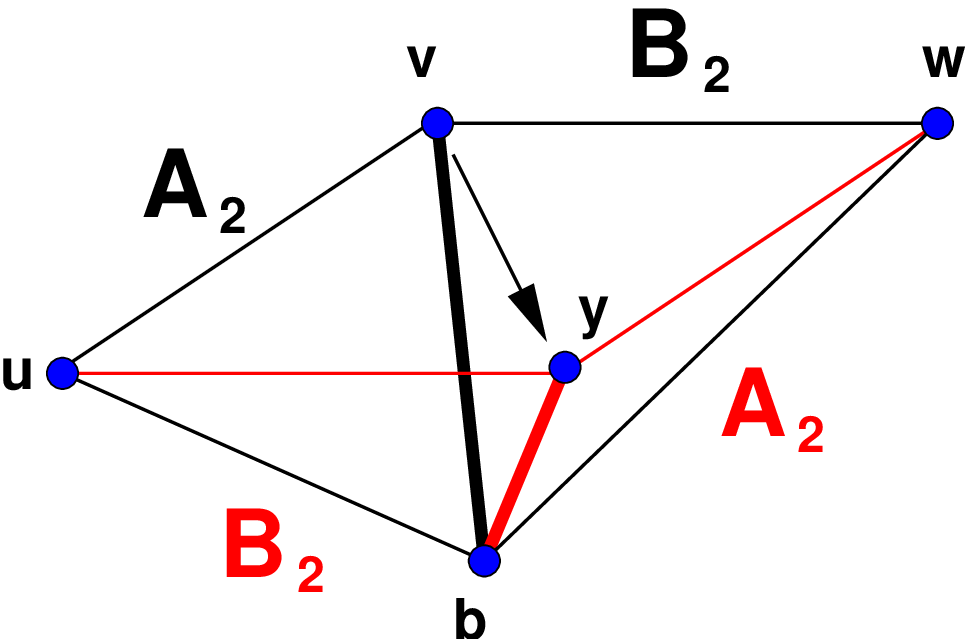}}}
\quad  A_2(b,u,v)B_2(b,v,w)&=&B_2(b,u,y)A_2(b,y,w), 
\quad y={u w+b\over v} \label{mutwo}
\end{eqnarray}
where we have represented in thick black (resp. red) line the initial (resp. mutated)
boundary edge. 
\end{lemma}

In the above equations, the transformations $b\to x$ (resp. $v\to y$)
are precisely the two types of forward mutations of the $A_2$ T-system cluster algebra, 
obtained by applying the first (resp. second) line of \eqref{tsysa2}.
We may now turn to the proof of Theorem \ref{genA2}.

\begin{proof}
The formula is proved by induction under mutation. We start from the basic staircase
solution \eqref{frisa2}, which may be put in the form \eqref{gena2}, upon substituting 
$A=A_1A_2$ and $B=B_1B_2$, and noting that
$k_1=1$, $j_1=j+k-1$, while $k_0=1$ and $j_0=j-k+1$. 

Starting from the expression \eqref{gena2} for the basic staircase boundary,
we may apply iteratively either of \eqref{mutone} or \eqref{mutwo} to get to any
other boundary (up to global translation), by simply substituting
products of pairs of triangle matrices into the expression \eqref{gena2}. 

We must however pay special attention to the extremal cases, namely
when the mutation acts on the edge just before the upper extremity as in 
Fig.\ref{fig:botop} (a), or just after the lower extremity as in Fig.\ref{fig:botop} (b), 
of the projection of $(1,j,k)$.

\begin{figure}
\centering
\includegraphics[width=12.cm]{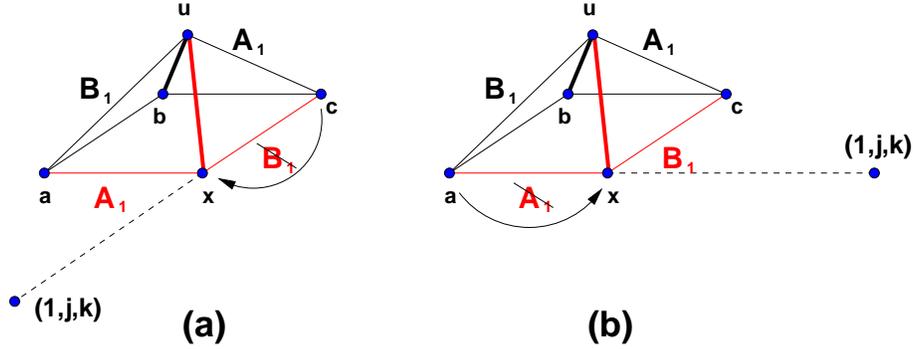}
\caption{\small The effect of a mutation next to the top (a) and bottom (b) of the
projection of $(1,j,k)$. The new mutated edge is represented in red. The corresponding
new lower vertex of the projection extremity is changed accordingly:
$c\to x$ (a) and $a\to x$ (b). The extremal matrices $B_1(x,c,u)$ (a) and $A_1(a,x,u)$ (b)
must be dropped from the expression for $M(i,j)$, as the corresponding triangles
lie outside of the new projection of $(1,j,k)$.}
\label{fig:botop}
\end{figure}

The prefactor $T_{1,j_1,k_1}$ in \eqref{gena2} corresponds indeed
to the bottom vertex of the upper extremity of the projection of $(1,j,k)$ onto the boundary.
Assuming as in Fig.\ref{fig:botop} (a) that the edge  just before the upper extremity 
of the projection of $(1,j,k)$ is of $e$ type, with value $b$ at the bottom vertex as 
in eq. \eqref{mutone}, the mutation
sends it to an $f$-type edge with bottom vertex value $x$, which becomes the new upper 
extremity of the projection of $(1,j,k)$, replacing $c$. 
Noting that $B_1(x,c,u)_{j,1}=\delta_{j,1} {x \over c}$, we see that the last multiplication
by $B_1(x,c,u)$ amounts to replacing the global prefactor $c$ by $x$, which is the
desired change of $T_{1,j_1,k_1}$. 
Analogously, when the mutation acts on an edge
of type $e$ next to the bottom extremity of the projection as in Fig.\ref{fig:botop} (b), 
with bottom vertex value $b$ as in \eqref{mutone},
it sends it to an edge of type $f$ with bottom value $x$, which becomes the new 
lower extremity of the projection of $(1,j,k)$, replacing $a$. As 
$A_1(a,x,u)_{1,j}=\delta_{j,1}$, we may drop the contribution of this first triangle, and we
recover \eqref{gena2}.
This completes the proof of the Theorem.
\end{proof}

The case $\al=2$ needs no extra work, due to the following symmetry:

\begin{lemma}
For any fixed boundary, with initial data of the form $\{T_{\al,j,k_{\al,j}}\}_{\al=1,2;j\in \Z}$
we have
$$T_{2,j,k}\left(\{T_{\al,j,k_{\al,j}}\}_{\al=1,2;j\in \Z} \right)
=T_{1,j,k}\left(\{T_{3-\al,j,k_{3-\al,j}}\}_{\al=1,2;j\in \Z}\right)$$
\end{lemma}
\begin{proof}
The transformation $\al\to r+1-\al$ is a symmetry of (\ref{Tsys}-\ref{bounT}).
\end{proof}

\begin{cor}\label{restgenA2}
The solution of the $A_2$ $T$-system for arbitrary boundary is for all $\al=1,2$
a Laurent polynomial of the initial data with non-negative integer coefficients.
\end{cor}

\section{The $A_r$ case}\label{mainsec}

\subsection{T-system}

We now consider the general $A_r$ T-system (\ref{Tsys}-\ref{bounT}) 
for $j,k\in \Z$, and $\al\in I_r$. 
As before, this splits again into two independent systems for $T_{\al,j,k}$
with fixed value of $\al+j+k$ modulo 2, which we fix to be $0$, 
without loss of generality. 

The indices $\{(\al,j,k)\}$ for $\al+j+k=0$ modulo $2$ 
run over $r$ consecutive horizontal layers $\al=1,2,...,r$
of the centered cubic lattice, each of which is a square lattice, the vertices of the next 
layer lying at the vertical of the centers of the faces of the previous one. 
For technical reasons, we will also represent the extra bottom and top layers 
$\al=0$ and $\al=r+1$, within which all values of $T$ are fixed to $1$, 
by the $A_r$ boundary condition \eqref{bounT}.

\subsection{Boundaries}
\label{stairsec}

\begin{figure}
\centering
\includegraphics[width=12.cm]{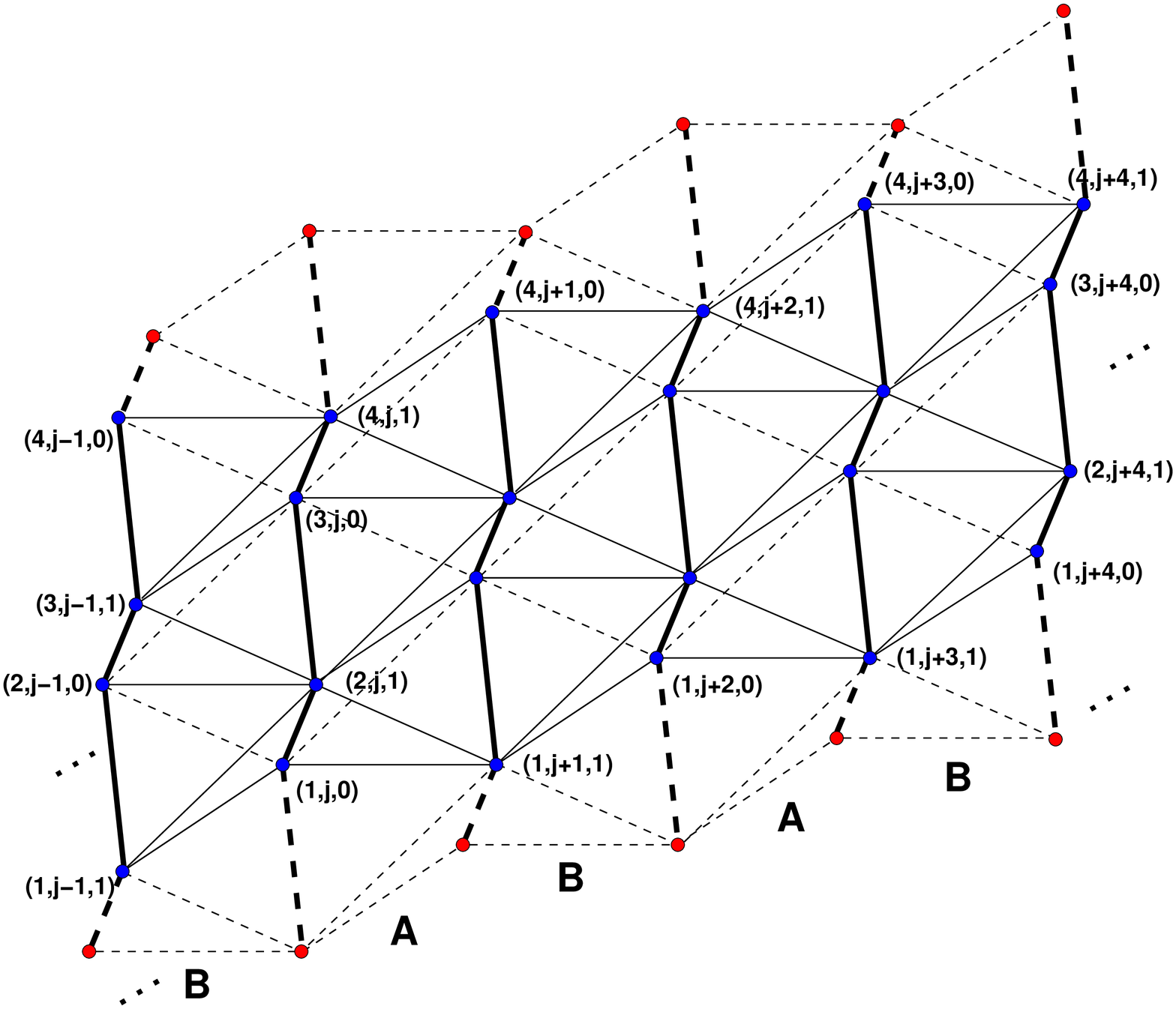}
\caption{\small The basic staircase boundary for the $A_r$ T-system, here for $r=4$, 
on which all vertex values $T_{\al,j,0}$ and $T_{\al,j-1,1}$ are specified (blue dots). 
We have added an
extra bottom and top layer where, according to the $A_r$ boundary,
 all vertex values are set to $1$ (red dots).}
\label{fig:stepar}
\end{figure}

As before we start with the basic boundary $\{T_{\al,j,0},T_{\al,j-1,1}\}_{\al\in I_r,j\in \Z}$ 
for $\al+j$ even. We may describe this boundary in 3D space as a succession of broken
lines at constant $j$ (represented in thick solid lines in Fig.\ref{fig:stepar}) of the form:
\begin{equation}\label{brok}
\ell_j=\{(\al,j,\epsilon_{\al,j})\}_{\al=1}^r, \quad \epsilon_{\al,j}=\al+j \, {\rm mod}\, 2
\end{equation}
We denote by A,B the vertical stacks of tetrahedrons depicted in Fig.\ref{fig:stepar}, 
respectively between $\ell_{2i-1}$ and $\ell_{2i}$
and $\ell_{2i}$ and $\ell_{2i+1}$ for all $i\in\Z$. 
Each tetrahedron has two opposite (thick) edges
of the form $(\al,j,\epsilon)-(\al+1,j,1-\epsilon)$ and $(\al,j+1,1-\epsilon)-(\al+1,j+1,\epsilon)$.

\subsection{Solution for the basic staircase boundary}

In Ref. \cite{DFK09a}, the system was solved for the basic staircase boundary
in two steps. First one eliminates $T_{\al,j,k}$ for all $\al>1$ as:
\begin{equation}\label{detalpha}
T_{\al,j,k}=\det_{1\leq a,b\leq \al} \, 
\left(T_{1,j-a+b,k+a+b-\al-1}\right), \quad \al\in I_r, \ j,k\in \Z
\end{equation}
which allows to concentrate on $T_{1,j,k}$.

The solution $T_{1,j,k}$ was then expressed in terms of paths on a rooted target
graph with time-dependent weights involving only the boundary values.
More precisely, $T_{1,j,k}$ was found to be equal to $T_{1,j+k,0}$ times
the partition function for paths starting from the root at time $j-k$ and ending at the
root at time $j+k$.
It is best expressed in terms of the $2r+2\times 2r+2$ transfer matrix, 
encoding the step weights $\by=y_1,y_2,...,y_{2r+1}$:
$$\cT(\by)=
\begin{pmatrix} 
0     & 1  & 0 & 0    & 0    & \cdots & 0 & 0 & 0 & 0 \\
y_1 & 0  & 1 & 0    & 0    & \cdots & 0 & 0 & 0 & 0 \\
0 & y_2  & 0  &  1 & 1 &  \cdots & 0 & 0 & 0 & 0 \\
0 & 0     & y_3  &  0    &  0   & \cdots & 0 & 0 & 0 & 0 \\
0 & 0     & y_4  &  0    &  0   & \cdots & 0 & 0 & 0 & 0 \\
\vdots &\vdots     &     &     &     & \ddots &  &  &  & \vdots \\
0 & 0     & 0    &  0    &  0   & \cdots & 0 & 1 & 1 & 0 \\
0 & 0     & 0    &  0    &  0   & \cdots & y_{2r-1} & 0           & 0         & 0 \\
0 & 0     & 0    &  0    &  0   & \cdots & y_{2r} & 0           & 0         & 1 \\
0 & 0     & 0    &  0    &  0   & \cdots & 0 & 0           & y_{2r+1}         & 0 
\end{pmatrix}
$$
Defining the time-dependent weights
$$ y_1(t)={T_{1,t,1}\over T_{1,t+1,0}} \quad 
y_{2\al+1}(t)={T_{\al+1,t-1,1}\, T_{\al,t+1,0}\over T_{\al+1,t,0}\, T_{\al,t,1}} \quad 
y_{2\al}(t)={T_{\al+1,t,1}\, T_{\al-1,t+1,0}\over T_{\al,t,0} \, T_{\al,t+1,1}}
\quad (\al=1,2,...,r)$$
and the transfer matrix for steps from time $t$ to $t+1$:
$$ \cT(t,t+1)=\cT(y_1(t),y_2(t),...,y_{2r+1}(t)) $$
we have:

\begin{thm}[\cite{DFK09a}]\label{thmar}
The solution of the $A_r$ $T$-system (\ref{Tsys}-\ref{bounT}) for the basic staircase boundary
reads for $\al=1$:
$$T_{1,j,k}=T_{1,j+k,0} \, \Big( \prod_{t=j-k}^{j+k-1} \cT(t,t+1) \Big)_{1,1} $$
\end{thm}

\subsection{Reduced transfer matrix}

As before we note that the two-step transfer matrix $\cT(\by,\by')=\cT(\by)\cT(\by')$
is again decomposable into two linear operators acting
on two complementary spaces of dimensions $r+1$,
corresponding respectively to components $(1,3,6,7,10,11,...,)$ and $(2,4,5,8,9,...)$. 
Explicitly, the operator acting on the first set of components reads:
$$\cP(\by,\by')= 
\begin{pmatrix} 
y_1'     & 1  & 0 & 0    & 0    & 0 & 0 & \cdots  \\
y_2y_1' & y_2+y_3'+y_4'  & 1 & 1    & 0    & 0 & 0 &  \\
0 & y_5y_4'  & y_5  &  y_5 & 0 &  0 & 0 &   \\
0 & y_6y_4'     & y_6  & y_6+y_7'+y_8'    & 1   & 1 & 0 &  \\
0 & 0     & 0  &  y_9y_8'   &  y_9   & y_9 & 0 &   \\
0 & 0     & 0  &  y_{10}y_8'   &  y_{10}   & y_{10}+y_{11}'+y_{12}' & 1 &    \\
\vdots &     &   &     &     & \ddots &  &   \\
\end{pmatrix}$$

The corresponding reduced two-step transfer matrix is $\cP(t,t+2)=\cP(\by(t),\by(t+1))$.
Theorem \ref{thmar} turns into:
\begin{equation}\label{cp}
T_{1,j,k}=T_{1,j+k,0}\, \Big( \prod_{i=0}^{k-1} \cP(j-k+2i,j-k+2i+2) \Big)_{1,1}
\end{equation}

\subsection{Gauge transformation and rhombus/triangle decomposition}

We may apply to $\cP$ any gauge transformation of the form
${\tilde \cP}(t,t+2)=L_t\cP(t,t+2)L_{t+2}^{-1}$ and such that 
$(L_t)_{1,j}=(L_t)_{j,1}=\delta_{j,1}$ without altering the result \eqref{cp}.

Here we choose:
$$L_t={\rm diag}(1,T_{1,t,0},{T_{2,t,1}T_{3,t,0}\over T_{3,t-1,1}},T_{3,t,0},\ldots,
T_{2\al-1,t,0},{T_{2\al,t,1}T_{2\al+1,t,0}\over T_{2\al+1,t-1,1}},\ldots)$$
One advantage of this choice is that ${\tilde \cP}(t,t+2)$ only depends
on values of $T_{\al,j,k}$ at times $j=t,t+1,t+2$. It is also justified {\it a posteriori}
by the decomposition formulas below.

We will now decompose the reduced two-step transfer matrix ${\tilde \cP}$
into a product 
of elementary matrices, defined as follows. 

\begin{defn}
We define the following $2\times 2$ elementary step matrices:
\begin{equation}\label{defHV}
H(a,b,x)=\begin{pmatrix} 1 & 0 \\ {x\over b} & {a\over b} \end{pmatrix}
\qquad  V(x,a,b)=
\begin{pmatrix} {a\over b} & {x\over b} \\ 0 & 1\end{pmatrix}
\end{equation}
\end{defn}
Note that these generalize the horizontal and vertical step matrices of Definition \ref{HVmat},
used in the $A_1$ case.

\begin{defn}
For any $2\times 2$ matrix $X$ and $\al\in \{1,2,...,r\}$ define the $r+1\times r+1$ matrices:
\begin{equation}\label{defX}
X_{\al,\al+1}=\left(\begin{array}{c|c|c} 
1_{\al-1} & 0_{\al-1\times 2} & 0_{\al-1\times r-\al}\\ \hline
0_{2\times \al-1} &  X & 0_{2\times r-\al} \\ \hline
0_{r-\al\times \al-1} & 0_{r-\al\times 2} & 1_{r-\al}
\end{array}\right)
\end{equation}
where $1_m$ denotes the $m\times m$ identity matrix and $0_{m\times p}$
the $m\times p$ matrix with zero entries.
\end{defn}

This gives rise to matrices $H_{\al,\al+1}(a,b,u)$ and $V_{\al,\al+1}(a,x,y)$ 
using for $X$ \eqref{defX} the matrices \eqref{defHV}. As an example, the 
triangle matrices $A_1,A_2,B_1,B_2$ introduced in the $A_2$ case \eqref{ABXY}
may be identified with:
\begin{equation*}
\begin{matrix} 
A_1(a,b,u)&=&H_{1,2}(a,b,u) \qquad\qquad & A_2(b,u,v)&=& V_{2,3}(b,u,v)\\
B_1(a,b,u)&=&V_{1,2}(1,a,b) \qquad\qquad & B_2(b,u,v)&=& H_{2,3}(u,v,1)
\end{matrix} 
\end{equation*}

We finally define:

\begin{defn}
We define the following $r+1\times r+1$ matrices:
\begin{eqnarray}
A_{\al}(a,x,y,u)&=&\left\{ \begin{matrix} 
H_{\al,\al+1}(x,y,u) & {\rm if}\ \al \ {\rm odd}\\
V_{\al,\al+1}(a,x,y) & {\rm if}\ \al \ {\rm even}\end{matrix} \right. \label{aimat} \\
B_{\al}(a,x,y,u)&=&\left\{ \begin{matrix} 
V_{\al,\al+1}(a,x,y) & {\rm if}\ \al \ {\rm odd}\\
H_{\al,\al+1}(x,y,u) & {\rm if}\ \al \ {\rm even}\end{matrix} \right. \label{bimat}
\end{eqnarray}
and for sequences $\ba=a_0,a_1,...,a_{r+1}$ and $\bb=b_0,b_1,...,b_{r+1}$,
we define:
\begin{eqnarray}\label{prodar}
A(\ba;\bb)&=&A_1(b_0,a_1,b_1,a_2)A_2(b_1,a_2,b_2,a_3)...A_r(b_{r-1},a_r,b_r,a_{r+1})
\nonumber \\
B(\ba;\bb)&=&B_1(b_0,a_1,b_1,a_2)B_2(b_1,a_2,b_2,a_3)...B_r(b_{r-1},a_r,b_r,a_{r+1})
\end{eqnarray}
\end{defn}

\begin{thm}
The following decomposition holds:
$$ {\tilde \cP}(t,t+2) =A(\ba(t),\ba(t+1))\,  B(\ba(t+1),\ba(t+2)) $$
where
$a_\al(t)=T_{\al,t,\epsilon_{\al,t}}$, with $\epsilon_{\al,t}$ as in \eqref{brok},
and in particular $a_0(t)=a_{r+1}(t)=1$, due to the $A_r$ boundary condition.
\end{thm}
\begin{proof}
By direct calculation.
\end{proof}

\begin{figure}
\centering
\includegraphics[width=15.cm]{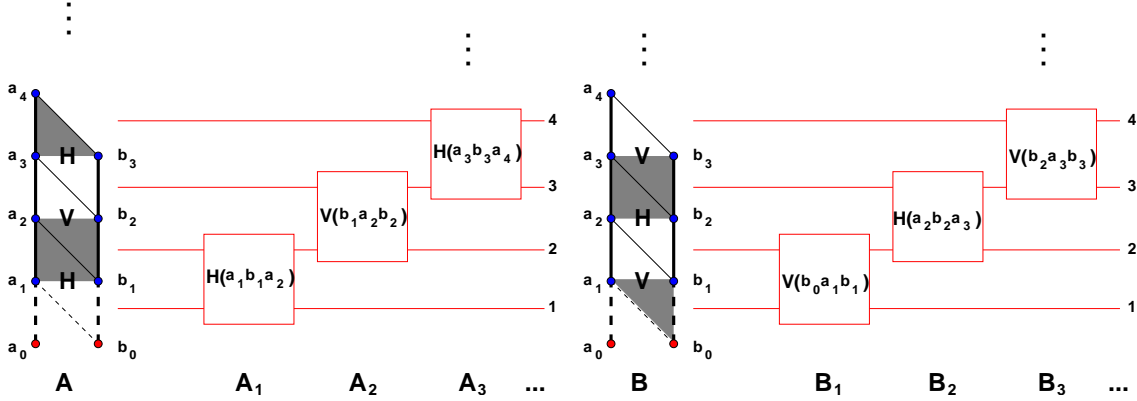}
\caption{\small The triangle decomposition of the stacks A,B of tetrahedrons on the
basic staircase boundary. The triangles are colored either white or gray, and are grouped 
into pairs of distinct colors sharing an edge (rhombi).
The H (resp. V) rhombi have the gray triangle pointing up (resp. down),
and correspond to products with the
matrices $H$ (resp. $V$) in the indicated manner. 
The three parameters are the values of $T_{\al,j,k}$
at the vertices of the gray triangle.}
\label{fig:abdecomp}
\end{figure}

In analogy with the frise result, we may now interpret pictorially the decomposition
result above. We interpret the matrices $A(\ba(j-1),\ba(j))$ and $B(\ba(j),\ba(j+1))$
as corresponding respectively to the stacks A,B of tetrahedrons in Fig.\ref{fig:stepar}.
Each tetrahedron may be decomposed in two ways into a pair of triangles 
(the thin solid lines in Fig.\ref{fig:stepar} correspond to one particular choice for each
tetrahedron). 

Let us now pick the canonical choice indicated in Fig.\ref{fig:abdecomp}
for both types A, B of stacks, namely that all squares are divided 
along their second diagonal. Then each term $A_i$ or $B_i$
in the products \eqref{prodar}
may be associated to a pair of triangles sharing an horizontal edge (which
we call a rhombus, by a slight abuse). 

As illustrated in Fig.\ref{fig:abdecomp}, we color the triangles in white or gray
in such a way that rhombi are made of one triangle of each color.
This gives rise to two types of rhombi H (resp. V), depending on whether the gray triangle
is on top (resp. bottom), that are interchanged when going from stacks of type A to 
type B.
In the rhombi, gray triangles play a particular role: their vertex values
carry the 3 parameters of 
respectively the matrices $H$ or $V$ in the product definition of
$A_i$ or $B_i$ (\ref{aimat}-\ref{bimat}). 
Taking the product over the rhombi from bottom to top yields the formulas \eqref{prodar}.

Each stack contains $r-1$ tetrahedrons, hence there are $2^{r-1}$ possible
triangle/rhombus decompositions of each stack of type A (resp. B), 
each made of $r$ rhombi. By their definition (\ref{defX}-\ref{aimat}-\ref{bimat}), 
the operators $A_i$ and $A_j$ (resp. $B_i$ and $B_j$) commute 
as soon as $|i-j|>1$. So a given rhombus decomposition carries the information of
whether $A_i$ multiplies $A_{i+1}$ from the left (like in Fig.\ref{fig:abdecomp}) 
or from the right (for the other choice of diagonal in the $i$-th tetrahedron from the bottom),
and of what the three arguments of the $H$ or $V$ factors are 
(via the boundary values at the vertices of the gray triangle). 

The $2^{r-1}$ {\it a priori} distinct matrix products corresponding to these rhombus
decompositions turn out to be identical. This is a consequence of the following 
local commutation relations (for odd $i$):
\begin{eqnarray*}
A_{i-1}(x,a,b,a')A_{i}(b,a',b',y)&=&A_{i}(a,a',b',y)A_{i-1}(x,a,b,b')\\
B_{i-1}(x,a,b,a')B_{i}(b,a',b',y)&=&B_{i}(a,a',b',y)B_{i-1}(x,a,b,b')
\end{eqnarray*}
easily derived from the following:

\begin{lemma}
\begin{eqnarray}
&& \quad \raisebox{-1.8cm}{\hbox{\epsfxsize=2cm \epsfbox{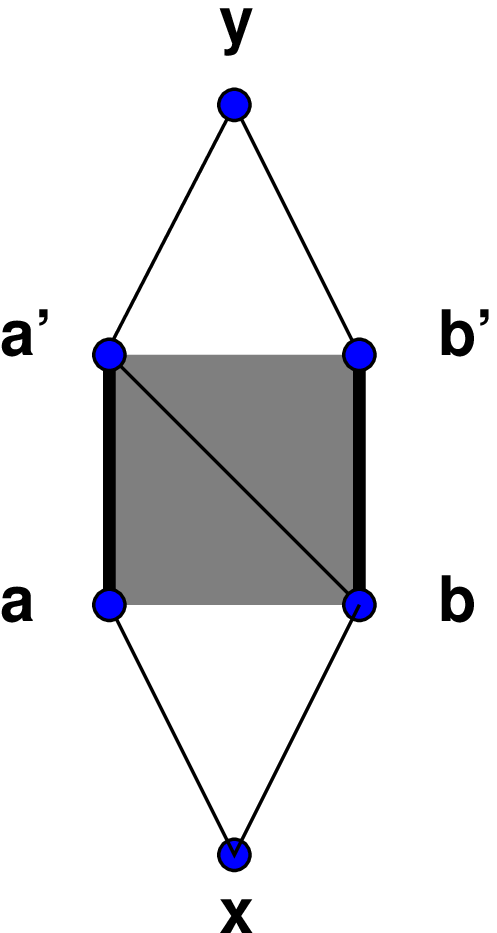}}}
\ H_{i-1,i}(a,b,a')\, V_{i,i+1}(b,a',b')=V_{i,i+1}(a,a',b')\, H_{i-1,i}(a,b,b')
\  \raisebox{-1.8cm}{\hbox{\epsfxsize=2cm \epsfbox{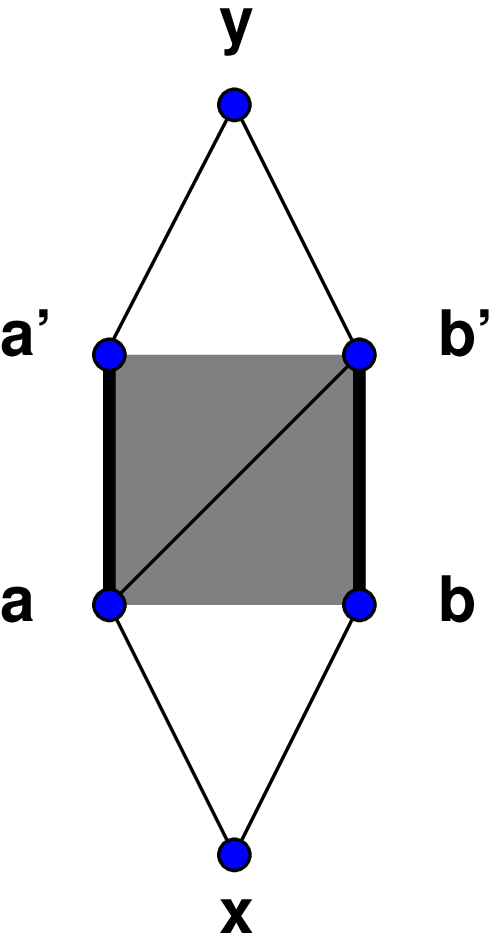}}}\label{aaid}\\
&& \quad \raisebox{-1.8cm}{\hbox{\epsfxsize=2cm \epsfbox{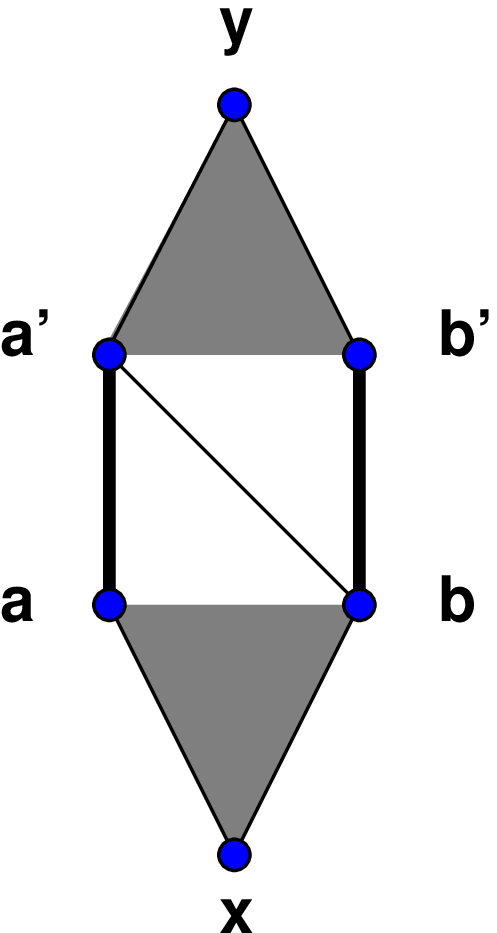}}}
\ V_{i-1,i}(x,a,b)\, H_{i,i+1}(a',b',y)=H_{i,i+1}(a',b',y)\, V_{i-1,i}(x,a,b)
\ \raisebox{-1.8cm}{\hbox{\epsfxsize=2cm \epsfbox{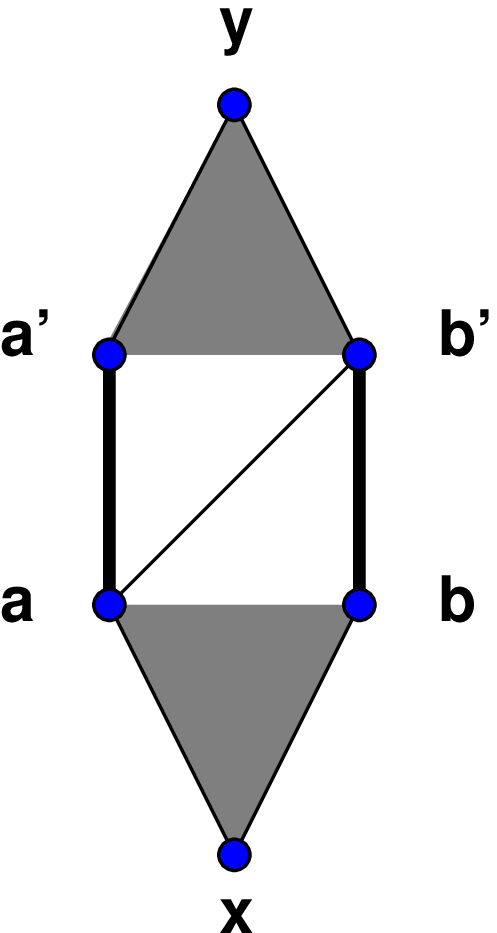}}}\label{bbid}
\end{eqnarray}
\end{lemma}
\begin{proof}
By direct computation.
\end{proof}

Theorem \ref{thmar} now turns into:
\begin{equation}\label{abt}
T_{1,j,k}=T_{1,j+k,0} \, \Big( \prod_{i=0}^{k-1} A(\ba(j-k+2i),\ba(j-k+2i+1))
B(\ba(j-k+2i+1),\ba(j-k+2i+2)) \Big)_{1,1}
\end{equation}
With the above-mentioned freedom of picking any of the $2^{r-1}$ $A$ and $B$ product
decompositions for each transfer matrix $A(\ba(t),\ba(t+1))$ and $B(\ba(t+1),\ba(t+2))$
in \eqref{abt}, we get $2^{k(r-1)}$ equivalent expressions for $T_{1,j,k}$.

\subsection{Mutations as rhombus exchange}
\begin{figure}
\centering
\includegraphics[width=11.cm]{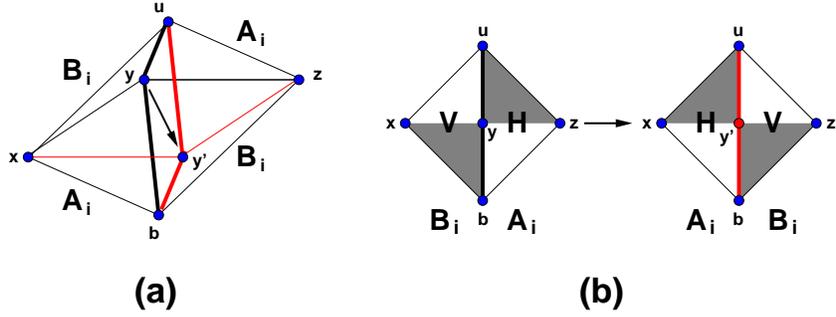}
\caption{\small A forward mutation of the boundary consists of a vertex value update
$y\to y'$, with $yy'=xz+au$. It takes the 4 triangles in the back (black lines) of the
octahedron (a) to the 4 triangles in the front (red lines).  
We have also represented (b) the corresponding $A_i$-$B_i$
matrix identity, equivalent to a rhombus exchange, in the case of odd $i$.}
\label{fig:mutar}
\end{figure}

Forward mutations of the boundary of the $A_r$ T-system
may take place whenever  five neighboring
vertices form the back of an octahedron, say in positions $(i,j,k-1)$,
$(i,j-1,k)$, $(i,j+1,k)$, $(i-1,j,k)$ and $(i+1,j,k)$, in which case the mutation 
$\mu_{i,j}$ replaces the back vertex $(i,j,k-1)$ with the front one $(i,j,k+1)$.
It also induces
the update of the boundary value $y$ of the back vertex
with $y'={bu+xz\over y}$ for the front one (see Fig.\ref{fig:mutar} (a), with $x=T_{i,j-1,k}$,
$y=T_{i,j,k-1}$, $z=T_{i,j+1,k}$, $b=T_{i-1,j,k}$, $u=T_{i+1,j,k}$, and $y'=T_{i,j,k+1}$).
Backward mutations just correspond to the inverse process with $y'\to y$.

Whenever a mutation is possible,
using the freedom to pick rhombus decompositions, we will see in next section that
we may always 
bring H and V rhombi in contact.
The mutation corresponds then
to interchanging the two rhombi: it is a forward mutation if V passes from the left to the
right of H (see Fig. \ref{fig:mutar}), backward otherwise.

To check that the mutation is faithfully represented by the
matrices, we use the following:

\begin{lemma}\label{switch}
\begin{equation}\label{identtwo}
\ \raisebox{-1.3cm}{\hbox{\epsfxsize=4.cm \epsfbox{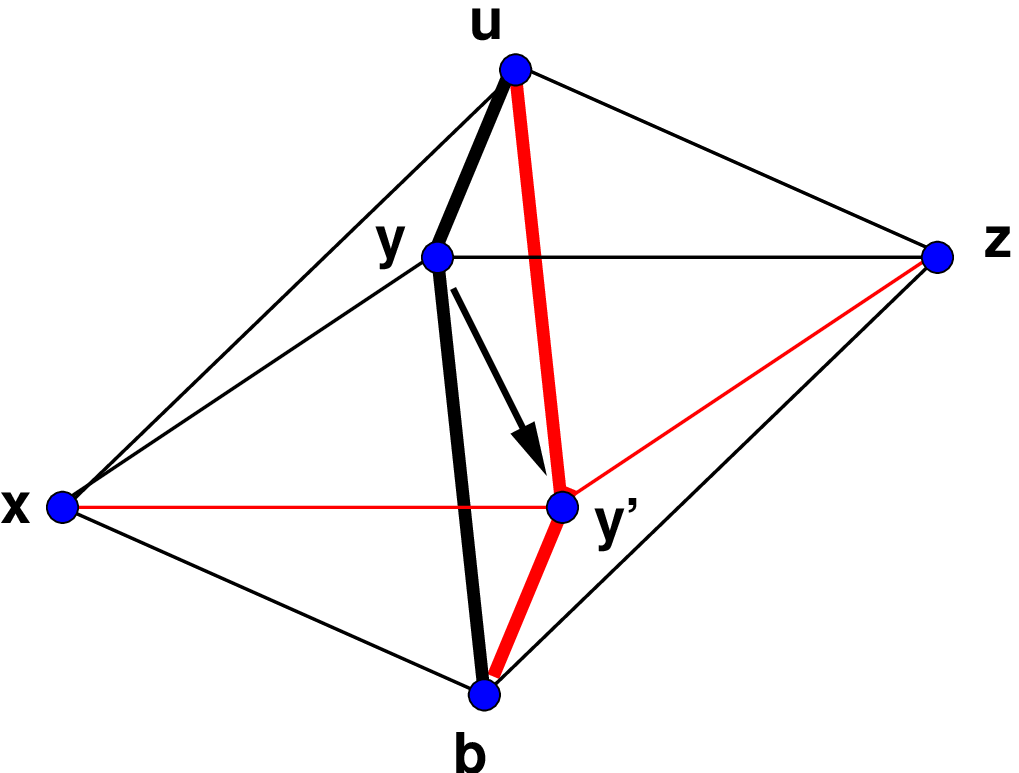}}}\
V(b,x,y)\, H(y,z,u)=H(x,y',u)\,V(b,y',z),\quad yy'=xz+bu
\end{equation}
\end{lemma}
\begin{proof} By direct calculation.
\end{proof}

The forward mutation $\mu_{i,j}$
is implemented
by simply switching the two rhombi in the decomposition of the corresponding
adjacent stacks. In the case $i$ odd, this corresponds to the identity depicted
in Fig.\ref{fig:mutar} (b):
$$ A_i(b,x,y,u)\, B_i(b,y,z,u)=B_i(b,x,y',u)\, A_i(b,y',z,u) \qquad yy'=xz+bu,$$
a direct consequence of Lemma \ref{switch}.
When $i$ is even, the mutation $\mu_{i,j}$ corresponds to the same identity 
read in the opposite direction
(but still corresponds to a V passing from the left to the right of an H).

We conclude that the forward mutations $\mu_{i,j}$ of the boundary may be iteratively
implemented on the formula \eqref{abt} by simply (i) picking the relevant 
rhombus decomposition among all equivalent ones namely bring a 
V to the left of an H with the vertex $(i,j,k_{i,j})$ of the boundary in the center (ii) switching 
the corresponding matrices $A_i$ and $B_i$ in the product transfer matrix. 
This will be made very precise in the next section.

\subsection{Arbitrary boundaries and the 6 Vertex model}

In this section, we give a bijection between boundary strips of the $A_r$ T-system and
configurations of the 6-Vertex (6V) model on infinite strips of width $r-1$. 
In the 6V picture, mutations
correspond to the reversal of all spins around square faces whose spins form an 
oriented loop.

\subsubsection{Local configurations and the six vertices}\label{conecsixv}

Recall that the 6V model is a statistical model defined on the square lattice, whose
configurations consist of an orientation (spin) on each edge, in such a way that 
the following ``ice rule" is satisfied: at each vertex of the lattice there are
exactly two entering and two outgoing edges.

The idea of the bijection is very simple. 
As we noted before, the boundary strips may be decomposed
into broken lines $\ell_j$ at time $j$, $j\in \Z$, as in \eqref{brok}, 
each made of a succession of $r-1$ edges of type $e$ or $f$. 
As discovered in the case of $A_2$, there are 6 possible transitions between an edge
of type $e$ or $f$ of the form
$e=(\al,j,k)-(\al+1,j,k+1)$ or $f=(\al,j,k)-(\al+1,j,k-1)$,
at time $j$, connecting vertices in the layers $\al$ and $\al+1$, and an edge 
of either type at time
$j+1$, also connecting vertices in the layers $\al$ and $\al+1$. 
The two corresponding
edges give rise to either tetrahedrons A,B if they are of opposite
types, and parallelograms C,D,E,F otherwise. 
Moreover, the 4 vertices included in the two edges
also define two horizontal edges respectively within the layers $\al$ and $\al+1$.
These may also be of only two types, say $e'$ and $f'$, of the form: 
$e'=(\al,j,k)-(\al,j+1,k+1)$ and $f'=(\al,j,k)-(\al,j+1,k-1)$.

\begin{figure}
\centering
\includegraphics[width=14.cm]{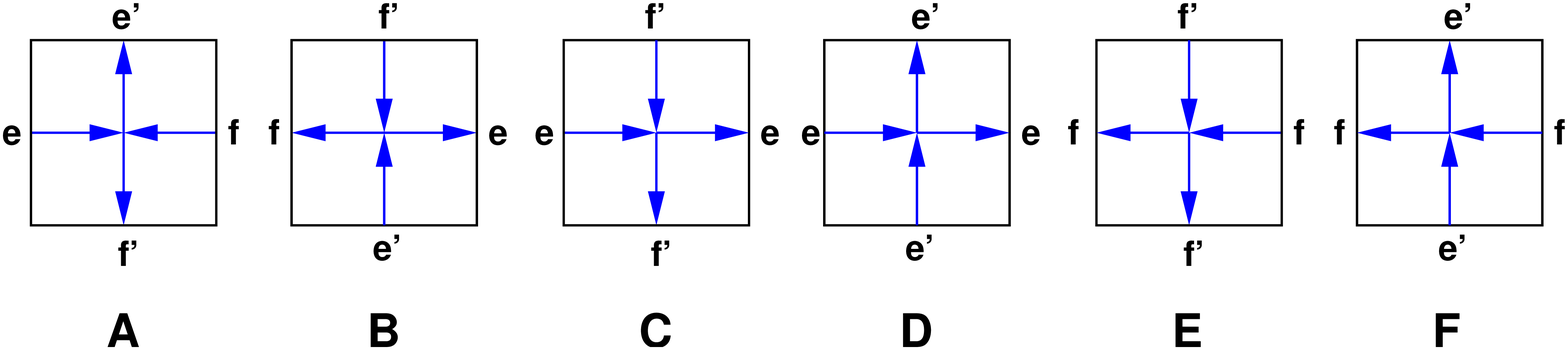}
\caption{\small The projections in $(\al,j)$ plane
of the six local configurations of the boundary of the $A_r$ T-system and the 
corresponding local arrow configurations of the dual 6V model. We have indicated
for each configuration the corresponding tetrahedron (A,B) or parallelogram (C,D,E,F).}
\label{fig:sixv}
\end{figure}

These 6 local configurations A,B,C,D,E may be represented as squares 
in projection onto a $(j,\al)$ plane, with vertical edges labeled $e$ or $f$
and horizontal edges labeled $e'$ or $f'$. In this projection,
a complete boundary is simply a strip of squares,
infinite in the $j$ direction, and of finite width in the $\al$ direction, with 
some compatible edge assignments. The corresponding configurations 
of the 6V model are on the dual of this strip, namely with a vertex in the middle of each
square.

To the 6 configurations of each square,
we
associate bijectively the 6 vertices of the 6V model as indicated in Fig.\ref{fig:sixv}, namely
we pick the horizontal edge orientation to be to the right (resp. left) 
for an $e$ type (resp. $f$ type) edge, while the vertical edge is oriented
upward (resp. downward) for an $e'$ type (resp. $f'$ type) edge.

Another direct way of connecting the initial boundary vertices to the 6V configuration
is to view the latter as coding the value of  $k_{\al,j}$ as the label of the face $(j,\al)$ in the
6V configuration as  follows. Note first that values of $k$ differ by $\pm 1$  in neighboring faces.
We impose the following ``Amp\`ere rule" that the value of $k$ on the right of an arrow is smaller than
that on the left. This fixes all the values of $k_{\al,j}$ up to a global translation.

\subsubsection{The basic staircase boundary as a 6V configuration}

\begin{figure}
\centering
\includegraphics[width=9.cm]{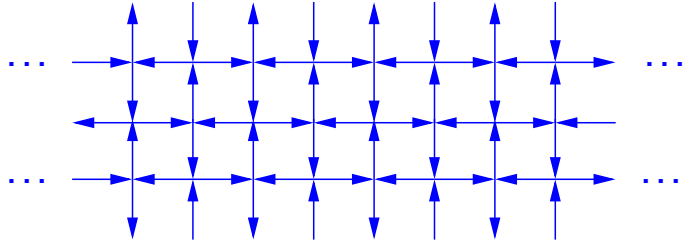}
\caption{\small The 6V configuration corresponding to the basic staircase configuration
of Fig.\ref{fig:stepar} for $r=4$.}
\label{fig:basicar}
\end{figure}

We now consider the basic staircase boundary introduced in Section \ref{stairsec}.
It is an infinite horizontal strip of square lattice with height $r-1$, with an alternance
of $...,e,f,e,f,...$ vertical edges and $...,e',f',e',f'...$ horizontal edges, namely a checkerboard
of configurations A and B of tetrahedrons, with vertices such that $k_{\al,j}\in\{0,1\}$. 
The dual 6V configuration
is simply made of arrows that alternate in all directions (it is a fully antiferromagnetic
groundstate configuration). In particular the bottom and top horizontal rows of vertical edges alternate 
along the $j$ direction.

\subsubsection{Mutations as loop reversals}

As illustrated in Fig.\ref{fig:mutar}, a mutation of the boundary at a given vertex
amounts to transforming the 4 edges sharing this vertex, by simply permuting
the two vertical ones (i.e. belonging to the plane $j=$const.) 
and the two horizontal ones (i.e. belonging to the plane $\al=$const.). 
This is also true if $\al=1$ or $\al=r$, with the condition that the bottom (resp. top)
vertex which belongs to the layer $\al=0$ (resp. $\al=r+1$) has an attached value $1$.
This situation was already encountered in the case $r=1$ (see Lemma \ref{mutaA1})
and $r=2$ (see Lemma \ref{boundamut}).
 
Such a mutation can only take place if
the permuted edges are of respective following types, clockwise around the 
vertex from the top: $e,e',f,f'$ for a forward mutation, and $f,f',e,e'$ for a 
backward mutation, and the two are transformed into one-another.

\begin{figure}
\centering
\includegraphics[width=14.cm]{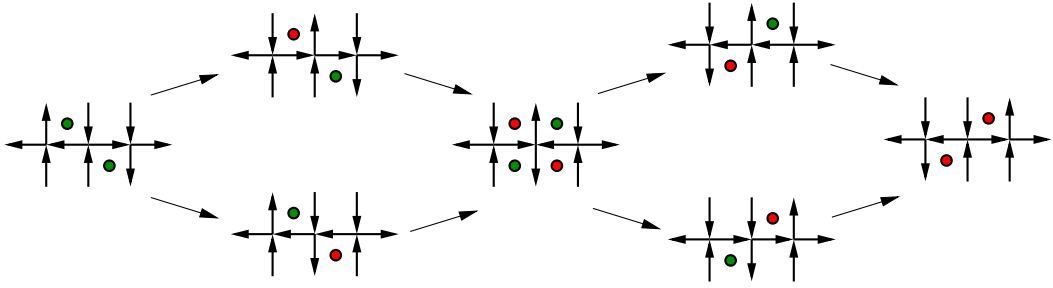}
\caption{\small A sequence of mutations applied on a length $2$ portion of the boundary
strip for $A_2$ in the 6V formulation.
The green (resp. red) dots mark the places where a forward (resp. backward) 
mutation can be applied. The spin of the top and bottom edge sequences are conserved, 
equal to $-1$ and $1$ here.}
\label{fig:loopreversal}
\end{figure}

 Once reformulated via the 6V bijection, and for $1<\al<r$,
 this simply means that a forward (resp. backward) mutation can take place
 only around faces whose adjacent edges form a clockwise (resp. counterclockwise) 
 oriented loop, and the mutation simply reverses the direction of all 4 arrows.
 When $\al=1$ (resp. $\al=r$), the ``loop" is open, i.e. is only formed of 3 edges
 with the bottom (resp. top) edge of the loop missing. 
 In the case of $A_1$, we are only left with open loops with both the top and bottom horizontal
 edges missing. In both cases, the corresponding mutation
 reverses the orientations of the 2 or 3 edges around the loop. By a slight abuse of language
 we still call ``loops" these open loops, and still refer to the edge reversal as ``loop reversal".

Starting from the basic configuration of Fig.\ref{fig:basicar}, we may therefore
generate all others by successive elementary loop reversals.
For illustration, we display in Fig.\ref{fig:loopreversal} a sequence of mutations
applied to a length $2$ portion of boundary in the $A_2$ case.
The central configuration corresponds to the basic staircase boundary.

All our 6V boundaries are generated by iterated loop reversals
on the 6V basic staircase. However, the basic staircase has the particular
property that its top and bottom boundary vertical edges alternate between 
up and down. Loop reversals involving these will create zones of successive
up or down edges, with the property that the total spin 
(i.e. the total number of up minus down edges) remains $0$. We deduce the following:

\begin{lemma}\label{boundasixv}
The boundaries of the $A_r$ T-system
correspond to the configurations of the 6V model on an infinite strip of width $r-1$, such that
there exist two times $j_{min}$ and $j_{max}\in \Z$ with the same parity such that:
\begin{itemize}
\item{(i)}
The portions
of 6V boundary for $j\leq j_{min}$
and for $j\geq j_{max}$ are identical to those corresponding to the basic staircase.
\item{(ii)}
The total spin for the portion for $j_{min}\leq j\leq j_{max}$ of 6V boundary is zero for both
top and bottom vertical edges.
\end{itemize}
\end{lemma}

Finally, the 6V boundaries must also carry the information of the initial data: there is one 
initial value per vertex of the original boundary, hence these may be represented
as face labels for the 6V
configuration (including the bottom and top faces with only 3 bounding edges).

\subsection{From general boundaries to networks}

We now wish to associate transfer matrices to the 6V configurations
described above. 

\subsubsection{6V and square/triangle decomposition}

\begin{figure}
\centering
\includegraphics[width=14.cm]{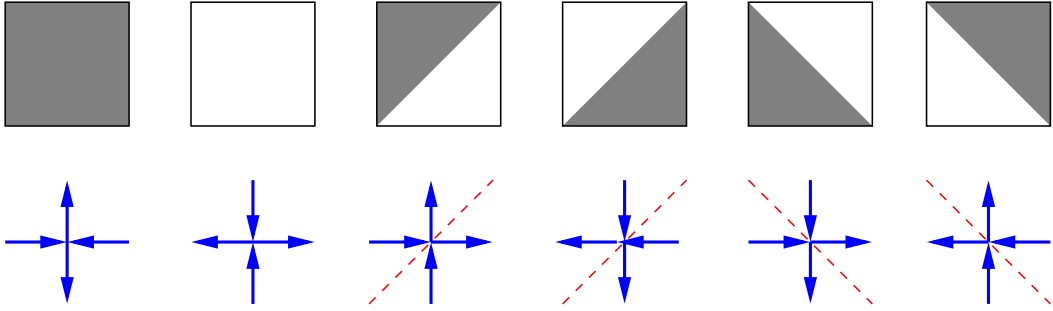}
\caption{\small The correspondence between 6V configurations 
and the square/triangle decomposition of the dual lattice. 
The diagonals that split squares into two triangles
are the diagonal reflection symmetry axes of the vertex configurations (dashed lines).
Vertices without diagonal reflection symmetry correspond to squares. The gray
(resp. white) triangles or squares
correspond to horizontal arrows pointing toward (resp. from) the central vertex,
and vertical arrows pointing from (resp. toward) the central vertex.}
\label{fig:sixvtriangle}
\end{figure}

To define the transfer matrices associated to 6V configurations, we first
associate to each 6V configuration a square/triangle decomposition of the dual
square lattice via the correspondence of Fig.\ref{fig:sixvtriangle}. 
Note that each edge
of the square lattice is adjacent to two polygons (triangle or square) of opposite colors.

\subsubsection{Slice transfer matrix}\label{slicetm}

\begin{figure}
\centering
\includegraphics[width=12.cm]{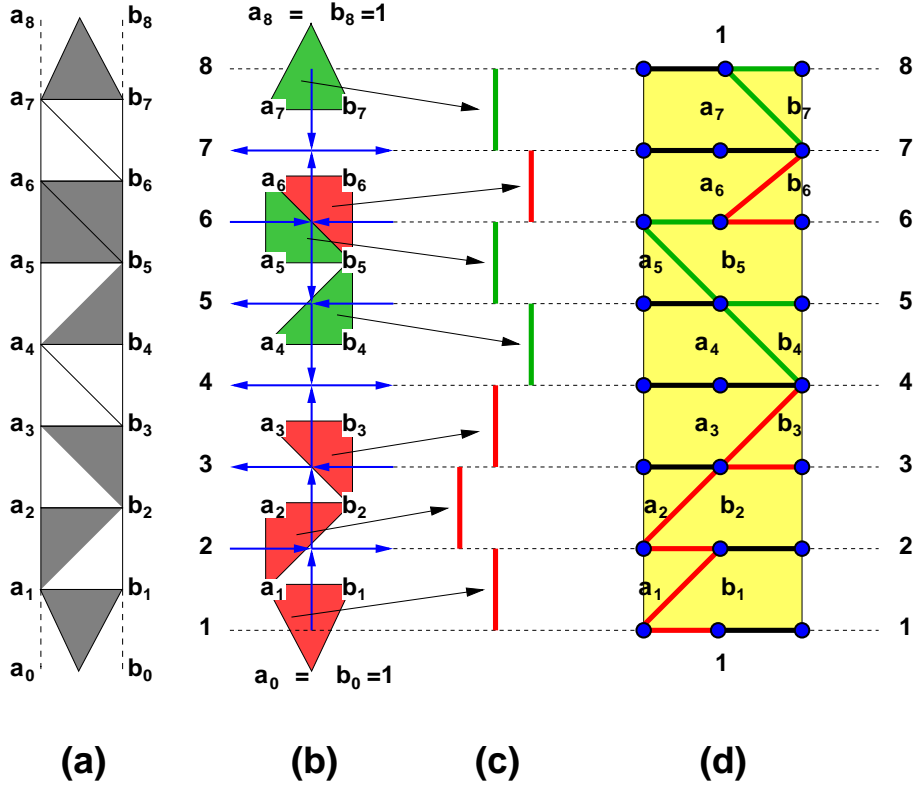}
\caption{\small A typical vertical slice in square/triangle decomposition (a)
and the associated 6V configuration (b) in the $A_7$ case. The squares in (a) are further
cut into two triangles by drawing the second diagonal. The configurations (a)-(b) 
determine the $H,V$ structure of the transfer matrix (c) where green (resp. red)
segments between lines $i,i+1$ correspond to $H_{i,i+1}$ (resp. $V_{i,i+1}$) 
matrices. The relative position of segments in consecutive horizontal slices
corresponds to the order of multiplication of the associated matrices.
We also indicate the network (d) corresponding to the complete matrix product.
The labels $a_i,b_i$ are now face labels of the network.}
\label{fig:column}
\end{figure}

Any given 6V boundary is decomposed into vertical slices
of unit width, made of a vertical sequence of $r-1$ vertices. To each such
slice we associate an $r+1\times r+1$ matrix $T$.
The matrix
$T$ is defined recursively
from the bottom to the top of the slice according to its 6V configurations as follows. 

For illustration and help, we have represented in Fig.\ref{fig:column} a typical
slice configuration (a) in square/triangle decomposition and (b) in 6V form.
We first consider the square/triangle decomposition of the 6V slice at hand, the
vertices of which carry the initial data, say $a_1,a_2,...,a_r$ from bottom to top
on the left and $b_1,b_2,...,b_r$ on the right.
We further split all the squares into pairs of triangles by drawing their second diagonal.
Note that horizontal edges are adjacent to two triangles of opposite colors. By a slight
abuse of language, we call again {\it rhombi} such pairs of triangles. With this definition, we
are left with two unpaired triangles with edges $a_r-b_r$ and $a_1-b_1$ respectively. 
This is repaired by adding the missing triangle of the opposite color on both top and 
bottom of the slice, with a third vertex carrying the value $a_{r+1}=b_{r+1}=1$ and
$a_0=b_0=1$. As they all carry the same value $1$,
it will be convenient to identify all the vertices of the bottom
layer as well as all those of the top one.

We label the lines supporting the horizontal edges of the 6V and the
two extra rows of vertices by $1,2,...,r+1$,
from bottom to top. 
We denote by $T_i$ the matrix corresponding to the partial
slice extending from lines $1$ to $i$. 

We start by defining $T_1=H_{1,2}(a_1,b_1,x)$ if the bottom triangle is white
and $T_1=V_{1,2}(1,a_1,b_1)$ if it is gray, where $x$ is the top vertex 
value of the rhombus containing the bottom triangle.
Let us consider the vertices $i$ and $i+1$, for $i=1,2,...,r-1$.
If the 6V vertical edge between line $i$ and $i+1$ points up (resp. down), 
then $T_{i+1}$ is obtained by multiplying $T_i$ by an operator $V_{i,i+1}$ 
(resp. $H_{i,i+1}$) of Eqs. (\ref{aaid}-\ref{bbid}). Moreover the order of
multiplication is from the left (resp. right) if the diagonal of the square is 
the first, connecting the SW and NE corners (resp. second, 
connecting NW to SE corners). This gives:
\begin{equation}\label{orderhv}
\epsfxsize=12cm \epsfbox{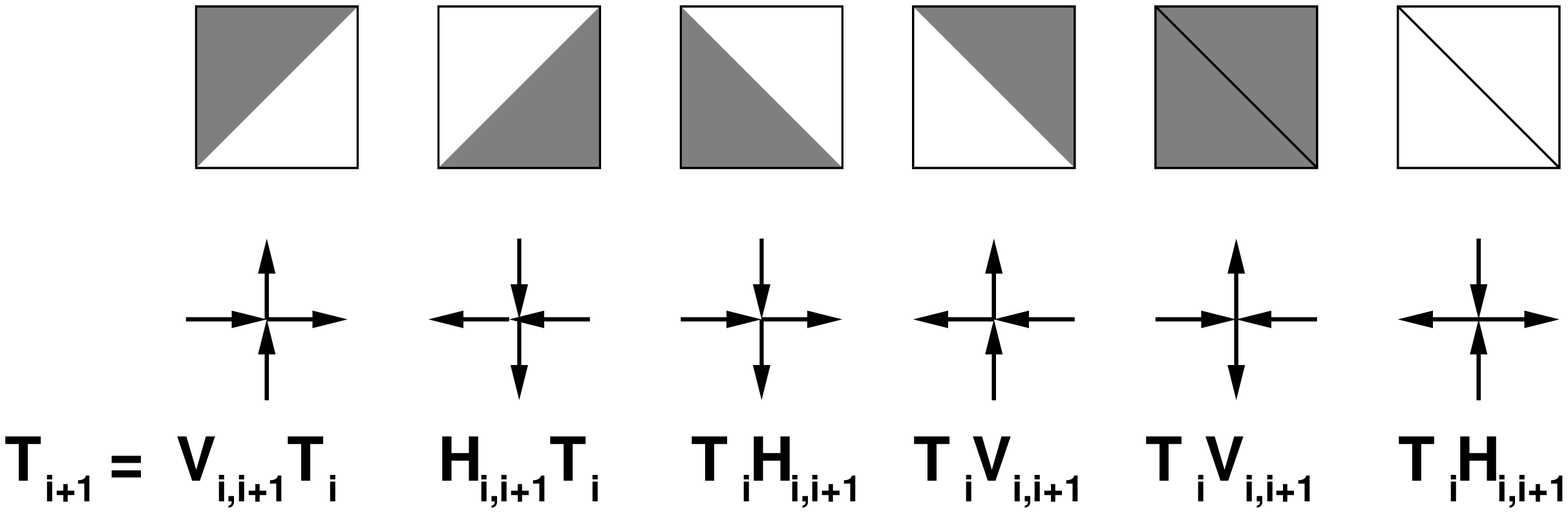}
\end{equation}
The three arguments of the $V_{i,i+1}$ (resp. $H_{i,i+1}$) operators are
$(x_{i-1},a_i,b_i)$ (resp. $(a_i,b_i,u_{i+1})$), where $x_{i-1}$ (resp. $u_{i+1}$)
is the third vertex label of the gray triangle with edge $a_i-b_i$.
This defines $T=T_{r}$ uniquely.

\begin{example}\label{netex}
Let us consider the case of Fig.\ref{fig:column} (a)-(c). 
The total matrix $T=T_7$ is obtained recursively
as follows:
\begin{equation*}\small \begin{matrix}
T_1=V_{1,2}(1,a_1,b_1)\quad \ \ &
T_2=V_{2,3}(a_1,a_2,b_2) \, T_1 &
T_3=T_2 \, V_{3,4}(b_2,a_3,b_3) &
T_4=T_3 \, H_{4,5}(a_4,b_4,b_5) \\
T_5=H_{5,6}(a_5,b_5,b_6) \, T_4 &
T_6=T_5 \, V_{6,7}(b_6,a_7,b_7) &
T_7=H_{7,8}(a_7,b_7,1) \, T_6  \end{matrix}
\end{equation*}
\end{example}

\subsubsection{Networks}

As pointed out earlier, the operators with respective indices $i,i+1$ and $j,j+1$ commute
as soon as $j>i+1$ or $j+1<i$. A useful way of representing the matrices 
$V_{i,i+1}$ and $H_{i,i+1}$ is as networks. 
We start with the representation:
\begin{eqnarray}\label{hvnet}
&&\epsfxsize=8cm \epsfbox{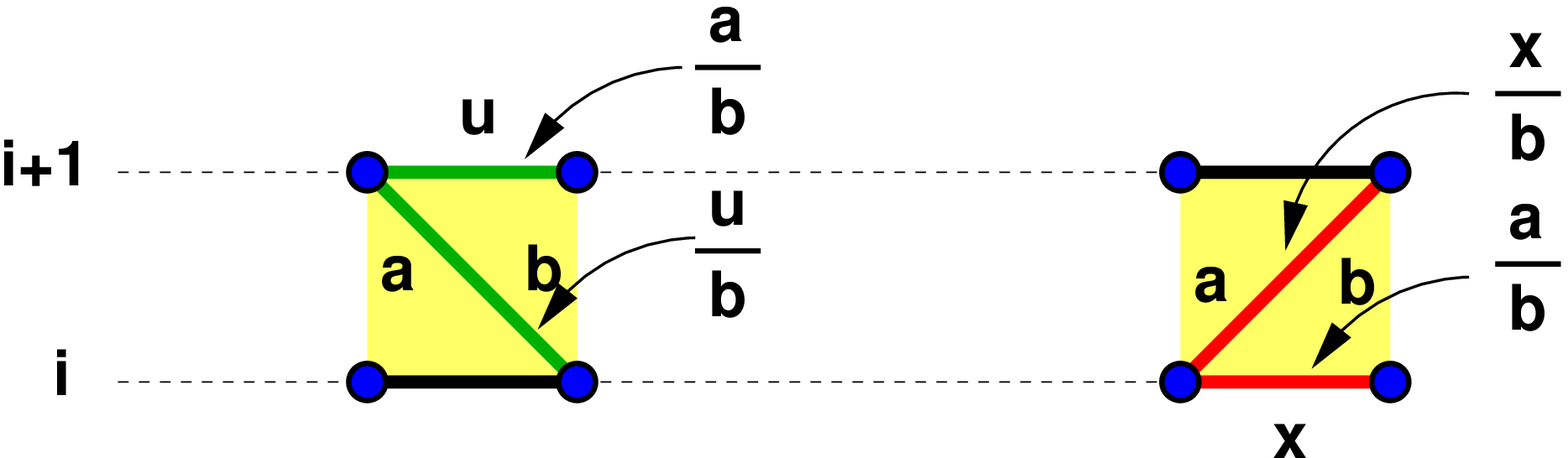} \nonumber \\
&& \qquad  \qquad \!\!\!  H_{i,i+1}(a,b,u) \qquad  \qquad 
V_{i,i+1}(x,a,b)
\end{eqnarray}
where the left vertices $i,i+1$ are connected to the right vertices $i,i+1$
via weighted edges $i\to j$ corresponding to the nonzero entries $(i,j)$ of the matrices $H$, $V$
of \eqref{defHV}. More precisely, a diagonal entry $1$ with indices $(j,j)$
is coded by  a thick black horizontal edge connecting the left vertex $j$ 
to the right vertex $j$. The other nontrivial entries $(j,k)$ are coded by thick green
(resp. red) edges for $H$ (resp. $V$), connecting the left vertex $j$ 
to the right vertex $k$. The actual values of these nontrivial entries
are coded by the labels of the faces of the network $x,a,b,u$ in \eqref{hvnet}.
Implicitly, the full representation of $H_{i,i+1}$ and
$V_{i,i+1}$ as $r+1\times r+1$ matrices
also involves thick horizontal black edges connecting the left vertices
$1,2,...,i-1,i+2,...,r+1$ to their right counterparts, corresponding to the 
block form of \eqref{defX},
and which we omitted in \eqref{hvnet}
for simplicity.

The multiplication of two matrices is naturally coded by the concatenation of their
networks, by identifying the right vertices of the network of the first matrix to the left
ones of that of the second. The result is in general a rectangle of height $r+1$,
and width possibly smaller than the number of matrices multiplied, as some horizontal
black lines may be freely removed to gain space.
We denote by $N(M)$ the network associated to the matrix $M$.
For illustration, we have represented in Fig.\ref{fig:column} (d)
the network $N(T)$ corresponding to the matrix $T$ of Example \ref{netex}.
As another example, the ``mutation" Lemma \ref{switch} reads in network language:
$$
\epsfxsize=6cm \epsfbox{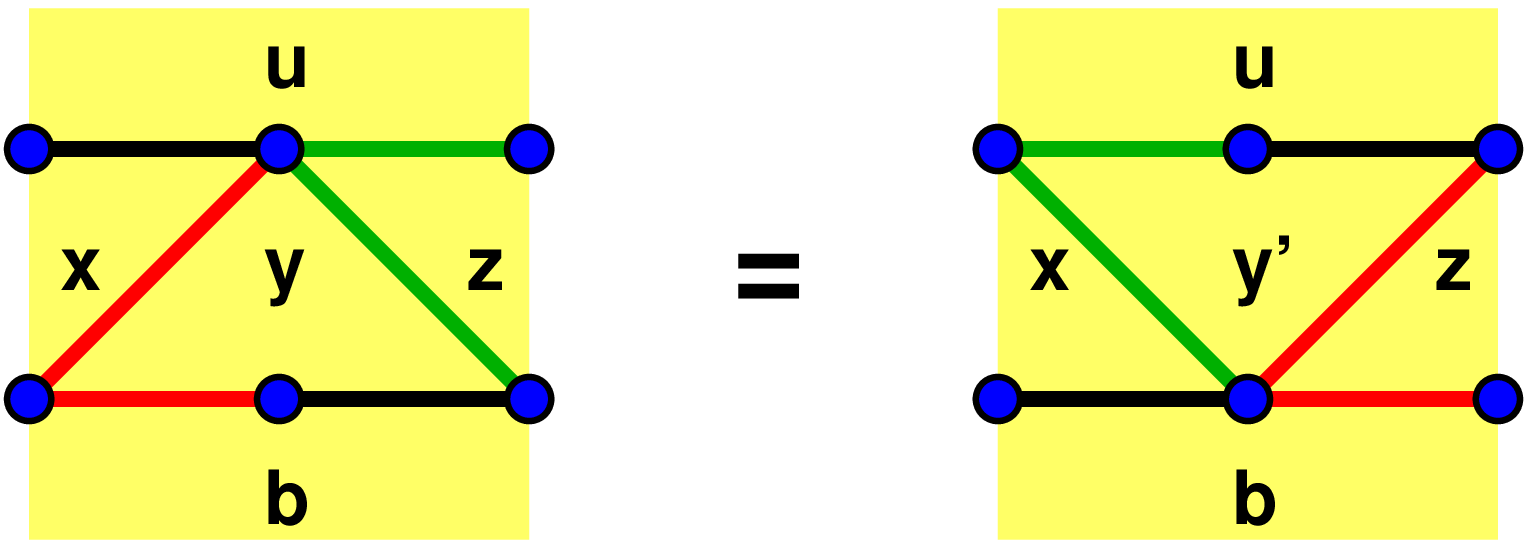}
$$

The networks provide us with natural weighted path models. 
Indeed, the matrix entry $T_{m,p}$ is interpreted as the partition function
for paths from left to right on the network of $T$, 
starting from the left vertex $m$ and ending 
at the right vertex $p$, weighted by the product of weights 
of the edges visited,
namely:
\begin{equation}
T_{m,p}=\sum_{{\rm paths}\ m\to p\atop
{\rm on}\ N(T)}  \prod_{{\rm edges}\ e \atop {\rm visited}} w(e)
\end{equation}

\subsection{General solution as path model on networks}

We now give the general expression for the solution $T_{\al,j,k}$ of the
$A_r$ T-system for arbitrary boundary conditions.

\subsubsection{The case of $T_{1,j,k}$}
As shown above, any boundary condition is coded by a configuration
of the 6V model on an infinite strip of width $r-1$ with the properties (i)-(ii)
of Lemma \ref{boundasixv}.

\begin{defn}
The projection of $(1,j,k)$ on the boundary 
is the portion of boundary between
the broken lines $\ell_{j_0}$ and $\ell_{j_1}$, 
containing the bottom vertices $(1,j_0,k_0)$ and $(1,j_1,k_1)$, 
repectively such that  $j_0-k_0=j-k$ and $j_1+k_1=j+k$ 
with $j_0$ maximal and $j_1$ minimal. 
\end{defn}

To this projection we naturally associate the corresponding truncated 
configuration $\cC$ of the 6V model on a rectangle of height $r-1$ 
between the planes $j=j_0$ and $j=j_1$, and with face labels 
given by the original vertex values
of the boundary.
Let $\cT_{j_0,j_1}(\cC)$ denote the product of vertical
slice transfer matrices from the slice $(j_0,j_0+1)$ to the slice $(j_{1}-1,j_1)$.
Then we have:

\begin{thm}\label{solonear}
The solution of the $A_r$ $T$-system with arbitrary fixed boundary reads for $\al=1$:
\begin{equation}\label{genear}
T_{1,j,k}=T_{1,j_1,k_1}\, \left(\cT_{j_0,j_1}(\cC)\right)_{1,1}
\end{equation}
in the above notations.
\end{thm}
\begin{proof}
This is proved by induction under mutation. First, \eqref{genear} is 
satisfied for $\cC=\cC_0$, corresponding to the truncated basic staircase boundary.
Indeed, we may rewrite \eqref{abt} as
\begin{eqnarray*} &&T_{1,j,k}=T_{1,j+k-1,1} \times \\ &&\left( \prod_{i=0}^{k-2} 
B(\ba(j-k+1+2i),\ba(j-k+2+2i)) A(\ba(j-k+2+2i),\ba(j-k+3+2i)) \right)_{1,1}
\end{eqnarray*}
by use of $A(\ba,\ba')_{1,j}=\delta_{j,1}$ and $B(\ba,\ba')_{j,1}=\delta_{j,1} a_1'/a_1$.
We see that this boils down to \eqref{genear}
with $k_0=k_1=1$, $j_1=j+k-1$,
and $j_0=j-k+1$. As explained above, any mutation $\mu$ may be implemented
by an elementary loop reversal on the 6V configurations, and only mutations 
within the projection of $(1,j,k)$ affect the value of $T_{1,j,k}$.

Let us assume \eqref{genear} holds for some 6V configuration with face labels $\cC$.
We wish to apply a mutation $\mu$, i.e. form the configuration $\mu(\cC)$,
identical to $\cC$ except for the reversed elementary loop and the updated face label.
Assume this is a forward mutation, i.e.
the corresponding loop is oriented clockwise in $\cC$, and assume for definiteness that
it occurs between vertical slices $s$ and $s+1$,
with respective slice
transfer matrices $T^{(s)}$ and $T^{(s+1)}$,
and between horizontal lines $i$ and $i+1$.
Due to the rules \eqref{orderhv}, examining the triangle decompositions
of the four adjacent vertices to the loop, we find the following
four possibilities:
$$
{\epsfxsize=12cm \epsfbox{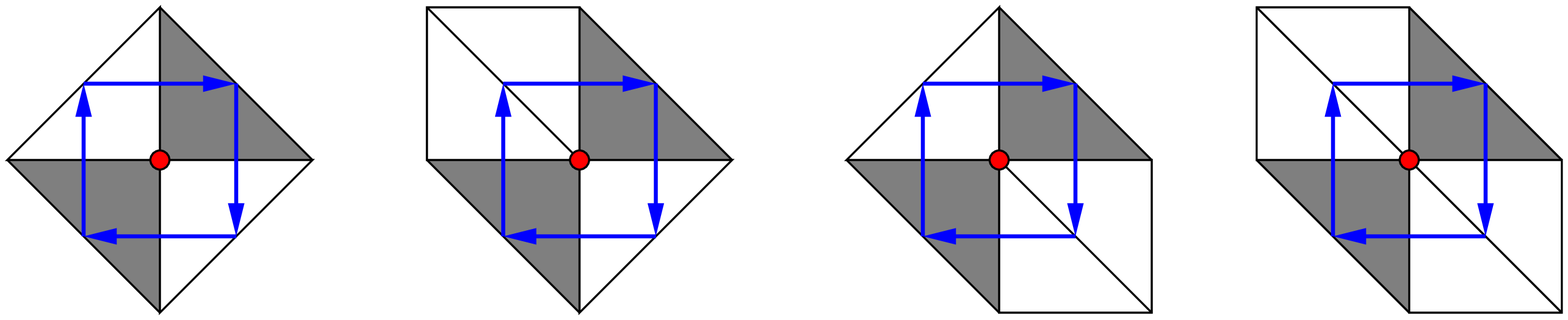}}
$$
Recalling that the choice of diagonal in the white squares is arbitrary,
we may bring all situations to the first one, by use of the identity \eqref{bbid}
within the relevant slice transfer matrices. Applying the mutation to this case
amounts to the transformation of Fig.\ref{fig:mutar} (b).
We now have the following four possibilities for the environment of the 
reversed loop:
$$
{\epsfxsize=12cm \epsfbox{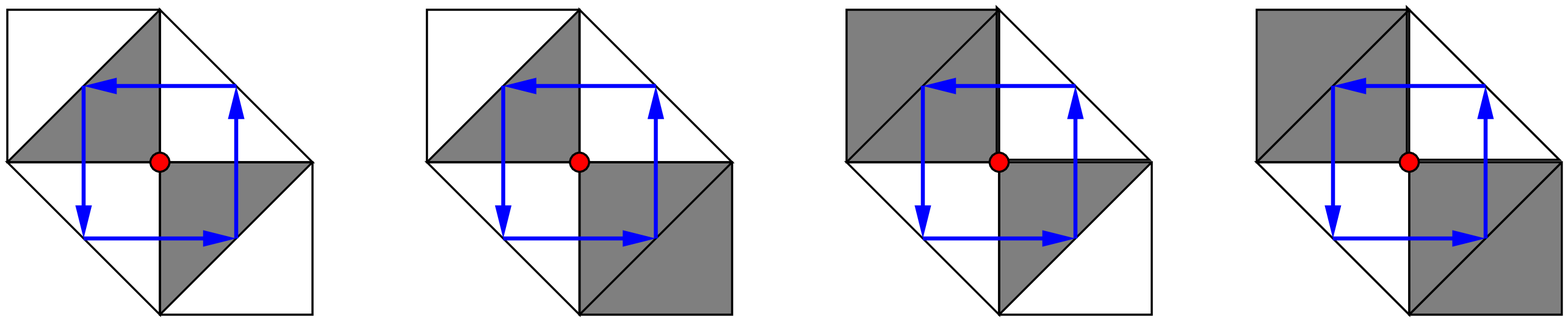}}
$$
To recover the correct slice transfer matrices, we must flip the diagonals
of the gray squares, by applying \eqref{aaid}. Taking the product
of transfer matrices over all slices, we finally get the ``mutated" matrix
$\cT_{j_0,j_1}(\mu(\cC))$.
We may repeat the same argument for backward mutations, while exchanging the roles
of white and gray triangles, and for mutations at vertices with $\al=1$ and $\al=r$ with
the obvious changes.

As before, we must however consider separately the case when 
the mutation occurs in the lower left
or lower right square face of the truncated 6V configuration. 
In these cases indeed, the projection of $(1,j,k)$ is modified and we must drop 
(resp. insert) a first or last slice transfer matrix when the mutation is forward 
(resp. backward). As before, dropping/inserting a last slice transfer matrix implements
the correct change of prefactor $T_{1,j_1,k_1}$ that makes it agree with 
the new projection.
In all cases, \eqref{genear} follows for $\mu(\cC)$. This completes the proof of
the Theorem.
\end{proof}

\subsubsection{Non-intersecting paths on networks and the case of $T_{\al,j,k}$}

Let $N_{j_0,j_1}(\cC)$ denote the network associated 
to the transfer matrix $\cT_{j_0,j_1}(\cC)$.
Then $T_{1,j,k}/T_{1,j_1,k_1}$ is interpreted as the
partition function for weighted paths on the network $N_{j_0,j_1}(\cC)$,
that start at the left vertex $1$ and end at the right vertex $1$, namely:
\begin{equation}\label{interpath}
{T_{1,j,k}\over T_{1,j_1,k_1}}=\sum_{{\rm paths}\ p:\ 1\to 1\atop
{\rm on}\ N_{j_0,j_1}(\cC)} \,  \prod_{{\rm edges}\ e \atop {\rm visited}} w(e)
\end{equation}
where $w(e)$ stands for the weight of the edge $e$.

Note that in the case of the basic staircase $\cC_0$,
this path interpretation is different from that of Ref. \cite{DFK09a}. 
Nevertheless, we may 
like in Ref. \cite{DFK09a} interpret the determinant identity \eqref{detalpha}
that relates $T_{\al,j,k}$ to $T_{1,j,k}$ in terms of non-intersecting paths,
now on networks.

Let us pick an arbitrary boundary $\cC$, and
rewrite the determinant formula  \eqref{detalpha} as:
\begin{equation}\label{newdet}
{T_{\al,j,k}\over \prod_{1\leq b\leq \al} T_{1,j_1(b),k_1(b)} }
=\det_{1\leq a,b\leq \al} \, 
\left({T_{1,j-a+b,k+a+b-\al-1}\over T_{1,j_1(b),k_1(b)}}\right)
\end{equation}
where we denote by $(1,j_1(b),k_1(b))$ the common rightmost bottom vertex of the 
projections onto the boundary $\cC$
of the vertices $(1,j-a+b,k+a+b-\al-1)$, $a=1,2,...,\al$. Similarly let us denote
by $(1,j_0(a),k_0(a))$ the common leftmost bottom vertex of the 
projections onto the boundary $\cC$ of the vertices 
$(1,j-a+b,k+a+b-\al-1)$, $b=1,2,...,\al$.

We denote by
$(j,x)$ the vertex of $N(\cC)$ that lies between the
slices $j-1,j$ and $j,j+1$ at height $x$, namely the common exit and entry
point $x$ of respectively the network for the slice $j-1,j$ and that for $j,j+1$.
For instance the left vertex $x$ of $N_{a,b}(\cC)$ has coordinates $(a,x)$,
while the right vertex of same height has coordinates $(b,x)$.

We finally have:

\begin{thm}\label{maintheo}
The solution $T_{\al,j,k}$ of the $A_r$ $T$-system for arbitrary boundary is equal
to $\prod_{1\leq b\leq \al} T_{1,j_1(b),k_1(b)}$ times the partition function for
families of $\al$
non-intersecting paths on the network $N(\cC)$, starting at 
the vertices $(j_0(a),1)$, $a=1,2,...,\al$ and ending at the vertices
$(j_1(b),1)$, $b=1,2,...,\al$.
\end{thm}
\begin{proof}
In view of \eqref{genear} and its reformulation in terms of paths on networks
\eqref{interpath},
we may interpret the $(a,b)$ term in the determinant \eqref{newdet}
as the partition function for paths from the bottom left to the bottom right vertex
on the network $N_{j_0(a),j_1(b)}(\cC)$. 
The Theorem follows from the application of the Lindstr\"om-Gessel-Viennot
theorem \cite{LGV1} \cite{LGV2}.
\end{proof}

\begin{cor}
The solution $T_{\al,j,k}$ of the $A_r$ $T$-system for arbitrary boundary 
condition $T_{\al,j,k_{\al,j}}=a_{\al,j}$ ($\al\in I_r,j\in\Z$) is a positive Laurent polynomial
of its initial data $(a_{\al,j})_{\al\in I_r;j\in\Z}$.
\end{cor}

\subsection{An example}

\begin{figure}
\centering
\includegraphics[width=14.cm]{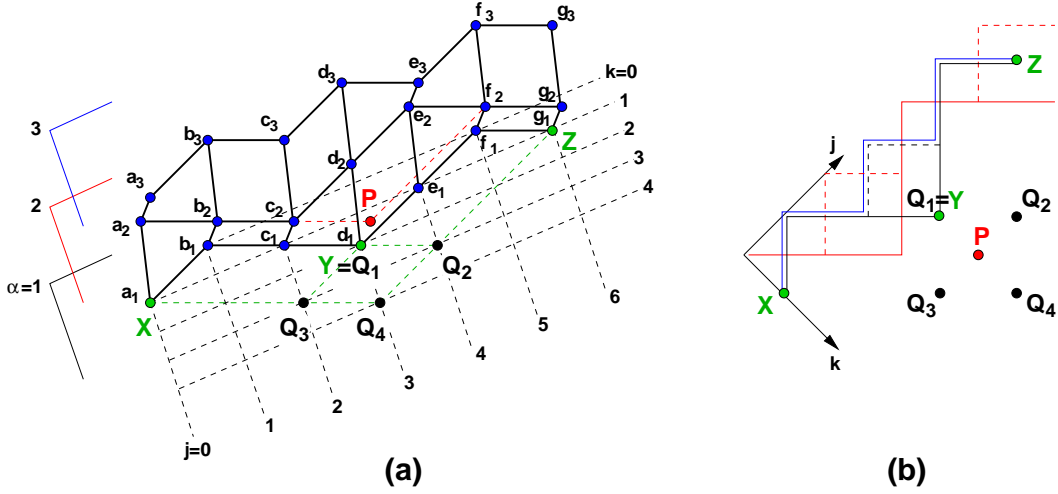}
\caption{\small A view in perspective
of a particular boundary (a) for the $A_3$ T-system (blue vertices
and thick solid black edges). We have indicated
the three horizontal layers $\al=1,2,3$, and the $j,k$ coordinates at $\al=1$.
We have represented the vertex $P=(2,3,3)$ and the four points of the $\al=1$ layer
$Q_1=(1,3,2),Q_2=(1,4,3),Q_3=(1,2,3),Q_4=(1,3,4)$ 
involved in the determinant formula
\eqref{deterQ}. The projections of the latter
onto the boundary are $X=(1,0,1)$, $Y=Q_1=(1,3,2)$ and $Z=(1,6,1)$. 
We also represent the same information (b) in projection onto the $(k,j)$ plane,
with layers $\al=1,2,3$ represented respectively in black, red, blue. We have
also indicated in dashed lines the mutations leading to this boundary from 
the basic staircase one.}
\label{fig:boundex}
\end{figure}

Let us consider the case $r=3$ and the boundary depicted in Fig.\ref{fig:boundex}, with
the $A_3$ boundary condition that $T_{0,j,k}=T_{4,j,k}=1$ for all $j,k$.
We wish to compute $T_{2,3,3}$ in terms of the indicated boundary values (blue vertices):
\begin{eqnarray*}
&&T_{1,0,1}=a_1 \ \ T_{1,1,0}=b_1 \ \ T_{1,2,1}=c_1 \ \ T_{1,3,2}=d_1 \ \ T_{1,4,1}=e_1 \ \ 
T_{1,5,0}=f_1 \ \ T_{1,6,1}=g_1 \\
&&T_{2,0,0}=a_2 \ \ T_{2,1,1}=b_2 \ \ T_{2,2,2}=c_2 \ \ T_{2,3,1}=d_2 \ \ T_{2,4,0}=e_2 \ \ 
T_{2,5,1}=f_2 \ \ T_{2,6,2}=g_2 \\
&&T_{3,0,1}=a_3 \ \ T_{3,1,0}=b_3 \ \ T_{3,2,1}=c_3 \ \ T_{3,3,0}=d_3 \ \ T_{3,4,1}=e_3 \ \ 
T_{3,5,0}=f_3 \ \ T_{3,6,1}=g_3
\end{eqnarray*}

We have the determinant formula
\begin{equation}\label{deterQ}
 T_{2,3,3}=\left\vert \begin{matrix} T_{1,3,2} & T_{1,4,3} \\
T_{1,2,3} & T_{1,3,4}\end{matrix}\right\vert=\left\vert \begin{matrix} T_{Q_1} & T_{Q_2} \\
T_{Q_3} & T_{Q_4}\end{matrix}\right\vert
\end{equation}
which involves the 4 points $Q_1,Q_2,Q_3,Q_4$
of the $\al=1$ plane.
In turn, the projections of the $Q_i$ onto the boundary read:
$X=(1,j_0(2),k_0(2))=(1,0,1)$, $Y=(1,j_0(1),k_0(1))=(1,j_1(1),k_1(1))=(1,3,2)$ and
$Z=(1,j_1(2),k_1(2))=(1,6,1)$.

\begin{figure}
\centering
\includegraphics[width=15.cm]{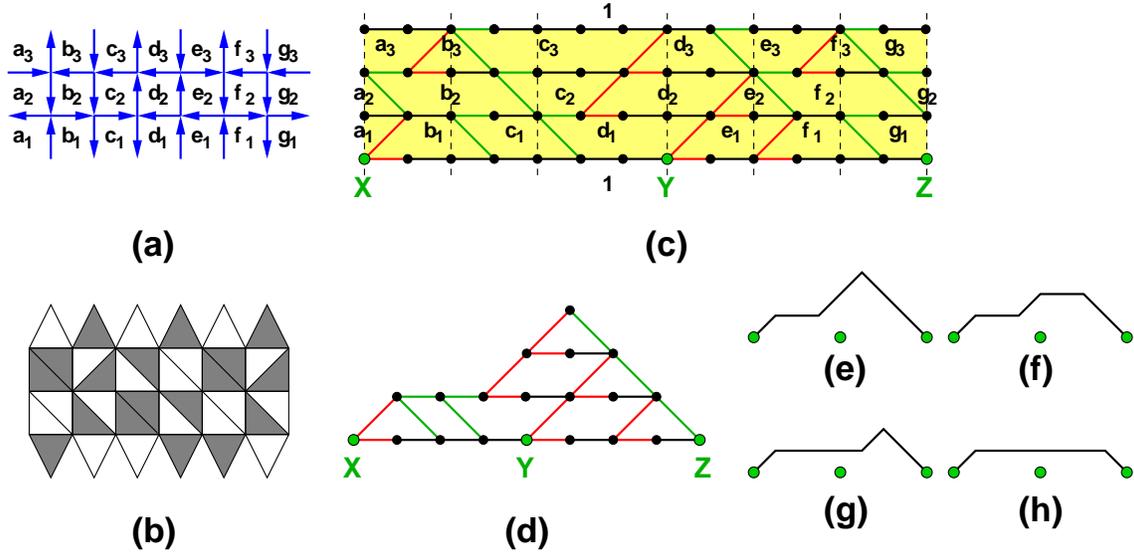}
\caption{\small The 6V version of the boundary of Fig.\ref{fig:boundex} (a)
with its face labels and the corresponding triangle decomposition of the dual (b).
The associated network (c) is represented with an indication of the successive slices
(vertical dashed lines) and the corresponding face labels. The truncation (d) of the
network is the only information needed to get the result. The partition function for
paths from $X$ to $Z$ that do not visit $Y$ is a weighted sum over the four paths 
(e),(f),(g),(h).}
\label{fig:exnetwork}
\end{figure}

Next, we construct the 6V configuration with face labels of Fig.\ref{fig:exnetwork} (a)
corresponding to the boundary. The simplest way is to start from the basic staircase
and apply to it the 3 mutations that bring it to the boundary under consideration
(see the dashed lines in Fig.\ref{fig:boundex} (b)). We also construct the triangle decomposition of Fig.\ref{fig:exnetwork} (b), by drawing the diagonal reflection
symmetry axes of the vertices and a second diagonal in the remaining squares.
This allows to build the network of Fig.\ref{fig:exnetwork} (c), by following the 
above procedure. In particular the edge-weights of the network are given
by the surrounding face labels according to the rules \eqref{hvnet}.

The Lindstr\"om-Gessel-Viennot theorem expresses
$T_{2,3,3}/(T_{1,3,2}T_{1,6,1})$ as the partition function for pairs of
non-intersecting paths on this
network that start at $X,Y$ and end at $Y,Z$. This is simply the partition function
for a single path from $X$ to $Z$ that does not visit $Y$. To compute the latter,
we may restrict ourselves to the truncation of the network depicted in 
Fig.\ref{fig:exnetwork} (d), where the edges keep their original weights from (c). 
The partition function is therefore a weighted sum over the four paths (e),(f),(g),(h),
and the final formula follows:
\begin{equation}
{T_{2,3,3}\over T_{1,3,2}T_{1,6,1}}={1\over d_3 g_1}+{c_3d_3\over d_2d_3g_1}+
{c_2e_1e_3\over d_1d_2e_2g_1}+{c_2f_2\over d_1e_2g_1}
\end{equation}

\section{Application: $Q$-system from the $T$-system solutions}

In this section, we compare the restriction of our $T$-system solution to the case of the $A_r$
$Q$-system whose solution was worked out in \cite{DFK3}. We find alternative path models
for describing the general solutions.

\subsection{$Q$ system}

The $Q$-system for $A_r$ is obtained from the $T$-system by ``forgetting"
about the variable $j$. This can be done in various ways, here we adopt the following:
we impose that the solution and the initial data be periodic, namely that $T_{\al,j+2,k}=T_{\al,j,k}$
for all $\al,j,k$.
In this case, the quantities
$R_{\al,k}=T_{\al,\al+k\, [2],k}$ are easily seen to satisfy the $A_r$ $Q$-system:
\begin{eqnarray}
R_{\al,k+1}R_{\al,k-1}&=&R_{\al,k}^2 +R_{\al+1,k}R_{\al-1,k}\nonumber \\
R_{0,k}&=&R_{r+1,k}=1\label{Qsys}
\end{eqnarray}

In \cite{DFK3}, the boundary data for \eqref{Qsys} were shown to be in bijection with
Motzkin paths of length $r-1$, $\bm=(m_1,m_2,...,m_r)$, and to take the form 
$x_\bm=\{R_{\al,m_\al},R_{\al,m_\al+1}\}_{\al=1}^r$. 
The explicit solution $R_{\al,k}$ for each such data was expressed as follows. 
First, for each Motzkin path $\bm$  one constructs an explicit rooted oriented graph $\Gamma_\bm$
with $2r+2$ vertices, and
with weighted edges encoded in a $2r+2\times 2r+2$ transfer matrix $T_\bm$. 
The root corresponds to the row and column index $1$ of the transfer matrix. The matrix $T_\bm$
is constructed in such a way that any oriented edge $i\to j$
of $\Gamma_\bm$ going away from the root (``ascent")
has a trivial weight $1$, while any oriented edge $i\to j$ 
pointing toward the root (``descent") has some-non-trivial weight
$t\times y_{i,j}(\bm)$, the latter being a Laurent monomial of the initial data.
The main result of \cite{DFK3} for $\al=1$ reads:
\begin{equation}\label{tobecomp}
R_{1,n+m_1}= R_{1,m_1} \left((I-T_\bm)^{-1}\right)_{1,1}\Big\vert_{t^n}
\end{equation}
where the notation $X\vert_{t^n}$ denotes the coefficient of $t^n$ in $X$. 
In other words, the quantity $R_{1,n+m_1}/R_{1,m_1}$ is the partition function for 
paths from and to the root on $\Gamma_\bm$, with exactly $n$ descents.

\begin{example}
For the case of $A_3$, the transfer matrix $T_{2,1,0}$ for the Motzkin path $(2,1,0)$ reads \cite{DFK3}:
\begin{equation}\label{tmatq}
T_{2,1,0}=\left( 
\begin{matrix}
0&1&0&0&0&0&0&0\\
t y_1&0&1&0&0&0&0&0\\
0&t y_2&0&1&1&0&0&0\\
0&0&ty_3&0&0&0&0&0\\
0&t y_{3,1}&ty_4&0&0&1&1&0\\
0&0&0&0&ty_5&0&0&0\\
0&t y_{4,1}&ty_{4,2}&0&ty_6&0&0&1\\
0&0&0&0&0&0&ty_7&0
\end{matrix} 
\right)
\end{equation}
where:
\begin{eqnarray}
&&y_1={R_{1,3}\over R_{1,2}},y_2={R_{2,2}^2\over R_{2,1}R_{1,2}R_{1,3}},
y_3={R_{1,2}R_{2,2}\over R_{2,1}R_{1,3}}, y_4={R_{1,2}^2R_{3,1}^2\over R_{3,0}R_{2,1}R_{2,2}R_{1,3}},
y_5={R_{2,1}R_{3,1}\over R_{3,0}R_{2,2}}\nonumber \\
&&y_6={R_{2,1}^2\over R_{3,0}R_{3,1}R_{2,2}},y_7={R_{3,0}\over R_{3,1}},
y_{4,2}={R_{1,2}^2\over R_{3,0}R_{2,2}R_{1,3}},y_{3,1}={R_{3,1}^2\over R_{3,0}R_{2,1}R_{1,3}},
y_{4,1}={1\over R_{3,0}R_{1,3}} \label{params}
\end{eqnarray}
\end{example}

The general solution for $R_{\al,n}$ was shown in \cite{DFK3} to be expressible as partition
function for families of $\al$ strongly non-intersecting weighted paths on $\Gamma_\bm$, via
a generalization of the Lindstr\"om-Gessel-Viennot Theorem.

\subsection{Network formulation from $T$-system}
As mentioned above, the $Q$-system restriction amounts to having a periodic solution
of the $T$-system, with period $2$ in the variable $j$.
For such a solution, boundaries are much simpler, as they are also 
$2$-periodic in the $j$ direction. However, constructing such a boundary by iterated mutations 
on the basic staircase (which has the desired periodicity)
involves repeating each mutation an infinite number of times
(with periodicity $2$ in the $j$ direction). So {\it stricto sensu} these boundaries are not
covered by our previous solution. But for each fixed value of $(\al,j,k)$
only a finite portion of the boundary is needed to express $T_{\al,j,k}$. So we may apply a finite
number of mutations for each case and get the correct result, and just formally complete the 
boundary by periodicity. We may therefore apply the results of Section \ref{mainsec} here.

In the periodic case,
the boundary values must read $T_{\al,0,k_{\al,0}}$ and $T_{\al,1,k_{\al,1}}$ for $j=0,1$ mod $2$ respectively.
Note that $k_{\al,1}-k_{\al,0}=\pm 1$.
Defining $m_\al={\rm Min}(k_{\al,0},k_{\al,1})$, $m_\al+1={\rm Max}(k_{\al,0},k_{\al,1})$, we
find that $\bm$ is the relevant Motzkin path for describing the $Q$-system boundary data.

According to our results, each periodic boundary may be viewed as a configuration
of the 6V model but now with a periodicity of $2$
in the direction of the strip. In other words, we have a configuration of the 6V model
on a cylinder of perimeter $2$ and height $r-1$. 
Moreover, as we started from a configuration with 
alternating vertical arrows (spin $0$ with two edges)
on both the
bottom and the top of the strip, any mutated configuration has one vertical edge pointing up 
and one pointing down on the upper and lower boundaries of the cylinder.
We have the following:

\begin{lemma}\label{motziv}
The 6V configurations
on a cylinder of perimeter $2$ and height $r-1$ with alternating edge orientations on top and bottom
are in bijection with Motzkin paths of length $r-1$, with $m_1=0$ or $1$.
\end{lemma}
\begin{proof}
The bijection goes as follows. We start from a Motzkin path $\bm=(m_1,m_2,...,m_r)$, with
$m_1=0$ or $1$. Then we have three possible situations for each of the $r-1$ steps of the Motzkin path:
$x_i=m_{i+1}-m_i=0,-1$, or $1$,
and two possible values for $y_i=i+m_i=0,1$ modulo $2$, for $i=1,2,...,r-1$. We have the following
dictionary between the six possible pairs $(x_i,y_i)$ and the six vertices of the 6V model.
$$
{\epsfxsize=14cm \epsfbox{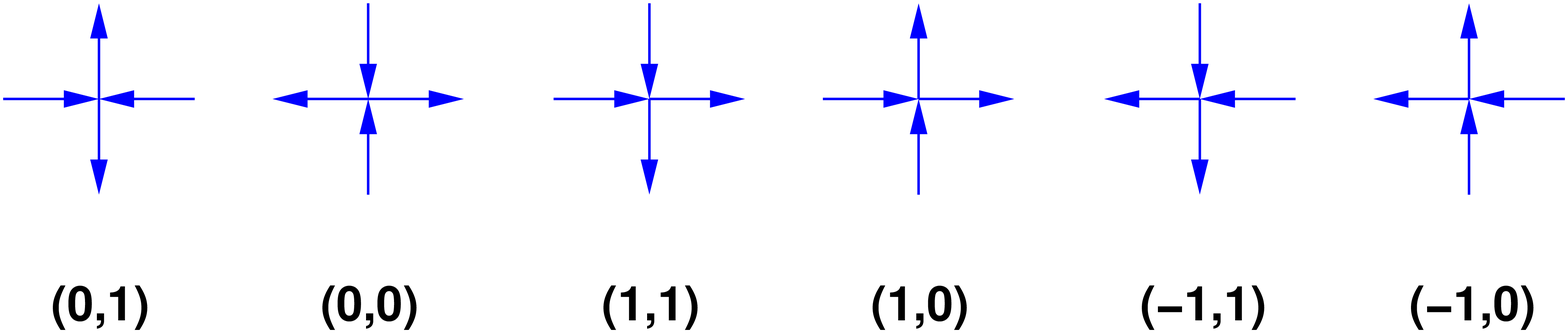}}
$$
More simply, $x_i$ is the ``algebraic sum" of the two horizontal edges ($0$ if they are opposed,
$1$ if they both point to the right, $-1$ if they both point to the left) and $y_i$ is determined only by the
bottom vertical edge ($1$ if it points down, $0$ if it points up).

Let us arrange the vertices associated to 
$\{(x_i,y_i)\}_{i=1}^{r-1}$ on a single column, from bottom to top. Note that the rules above make
the orientations of edges compatible, so we can identify the bottom vertical 
edge of the vertex $i+1$ with the top vertical edge of the vertex $i$. 
There is a unique way to complete this configuration 
into one on the $2\times r-1$ cylinder. 
We must indeed add a second column of vertices which are identical to the previous ones, 
up to reversal of all vertical arrows (the unique solution respecting the ice rule), 
in order to satisfy both the horizontal periodic boundary condition and the alternating one on top and bottom. 
Note that $m_1=0$ or $1$ determines whether the bottom left vertical edge points down or up.
\end{proof}

\begin{figure}
\centering
\includegraphics[width=10.cm]{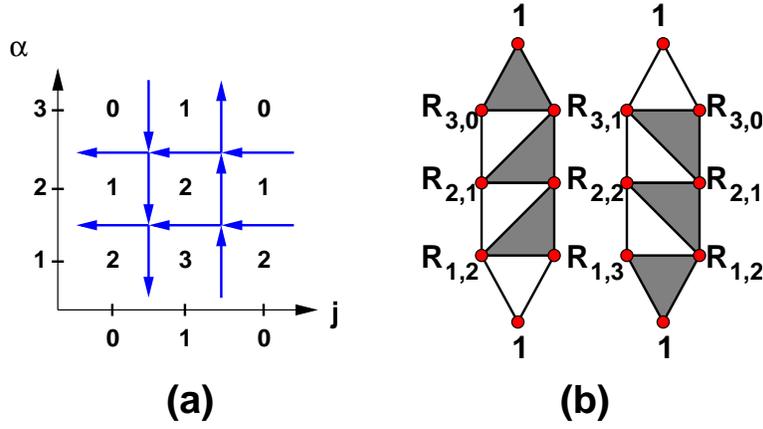}
\caption{\small The periodic 6V configuration (a) associated to the Motzkin path $(2,1,0)$ for $r=3$.
The corresponding boundary vertices $(\al,j,k)$ may be read off the configuration: the third
index $k$ is determined by
the Amp\`ere rule and by $m_1=2$, and is indicated in each face. We have also sketched the
associated rhombus/triangle decomposition (b) of the two corresponding slices, 
with the assigned boundary values.}
\label{fig:exqsys}
\end{figure}

\begin{example}
We consider the case $r=3$ and the Motzkin path $(2,1,0)$. We have $x_i=-1,y_i=1$ for $i=1,2$.
The 6V configuration is represented in Fig.\ref{fig:exqsys} (a). We have indicated the values of
$k_{\al,j}$ in the faces of the configuration, corresponding to $\al=1,2,3$ and $j=0,1$. 
These are not to be mistaken for the usual face labels of the 6V configuration, which are the assigned
boundary data $a_{\al,j}$ for the corresponding vertices: $T_{\al,j,k_{\al,j}}=a_{\al,j}$.
\end{example}

To each Motzkin path $\bm=(m_1,m_2,...,m_r)$, we may now associate
a configuration of the $6V$ model on a cylinder of perimeter 2 and height $r-1$ via Lemma \ref{motziv}
by using the shifted Motzkin path $\bm'=\bm-2(p,p,...,p)$, where $p=[{m_1\over 2}]$. The corresponding values
of $k_{\al,j}$ are uniquely determined by the Amp\`ere rule and by
$m_\al={\rm Min}(k_{\al,0},k_{\al,1})$, $m_\al+1={\rm Max}(k_{\al,0},k_{\al,1})$.
The solutions for $T_{1,j,k}$ of Theorem \ref{solonear} and for $T_{\al,j,k}$ of Theorem \ref{maintheo}
involve only two types of slice transfer matrices, due to the periodicity, namely
that for all slices of the form $[2a,2a+1]$ and that for $[2a-1,2a]$. 
These are coded by the two columns of the
6V configuration, which have the same horizontal edge orientations 
and opposite vertical ones. Using the construction of Sect. \ref{slicetm},
the slice transfer matrices may be constructed inductively as follows, directly from the Motzkin path $\bm$:

We start from $m_1$. If it is even, the first slice transfer matrix has $T_1=H_{1,2}$, if it is odd
it has $T_1=V_{1,2}$. Then, having constructed $T_i$, say with a last factor
$H_{i,i+1}$ (resp. $V_{i,i+1}$), we have three possibilities (according to the value of $x_i=m_{i+1}-m_i$): 
\begin{itemize}
\item{(i)} $x_i=0$: then $T_{i+1}=V_{i,i+1}T_i=T_iV_{i,i+1}$ (resp. 
$T_{i+1}=H_{i,i+1}T_i=T_iH_{i,i+1}$).
\item{(ii)} $x_i=1$: then $T_{i+1}=T_iH_{i,i+1}$ (resp. 
$T_{i+1}=V_{i,i+1}T_i$).
\item{(iii)} $x_i=-1$: then $T_{i+1}=H_{i,i+1}T_i$ (resp.
$T_{i+1}=T_iV_{i,i+1}$).
\end{itemize}
The arguments of the matrices are, as usual, the boundary data at the vertices of the gray triangles,
with the vertex indices of the form $(\al,j,k_{\al,j})$ where $k_{\al,j}$ are determined
by $\bm$. This gives a matrix $U_\bm=T_r$ for each Motzkin path $\bm$. The second slice is treated
analogously. 
Due to the fact that the its 6V configuration is identical to the first up
to reversal of all vertical arrow, it is easy to 
write the corresponding slice transfer matrix ${\tilde U}_\bm$
in terms of $U_\bm$. Indeed, the rhombus/triangle decomposition of the second slice is the reflection
of the first w.r.t. a vertical axis, and with all colors of triangles inverted. We therefore
have to interchange $H\leftrightarrow V$, and to reverse  the order of the factors. More precisely,
let $*$ be the involutive anti-automorphism ($(ab)^*=b^* a^*$, $(a^*)^*=a$) 
such that $V_{i,i+1}^*=H_{i,i+1}$, then we have ${\tilde U}_\bm=(U_\bm)^*$.

\begin{example}
In the case of the Motzkin path $(2,1,0)$ for $r=3$, we easily read the transfer matrices
on the rhombus/triangle decomposition of Fig.\ref{fig:exqsys} (b):
\begin{eqnarray}
U_{2,1,0}&=&H_{3,4}(R_{3,0},R_{3,1},1)H_{2,3}(R_{2,1},R_{2,2},R_{3,1}) 
H_{1,2}(R_{1,2},R_{1,3},R_{2,2})\nonumber \\
&=& \begin{pmatrix} 
1 & 0 & 0 & 0 \\
{R_{2,2}\over R_{1,3} } & {R_{1,2}\over R_{1,3} } & 0 & 0\\
{R_{3,1}\over R_{1,3} } & {R_{1,2}R_{3,1}\over R_{1,3}R_{2,2} }&{R_{2,1}\over R_{2,2} } & 0\\
{1\over R_{1,3}} & {R_{1,2}\over R_{1,3}R_{2,2} }&{R_{2,1}\over R_{2,2} R_{3,1}}&{R_{3,0}\over R_{3,1}}
\end{pmatrix}
\nonumber \\
{\tilde U}_{2,1,0}&=&V_{1,2}(1,R_{1,3},R_{1,2})
V_{2,3}(R_{1,2},R_{2,2},R_{2,1})V_{3,4}(R_{2,1},R_{3,1},R_{3,0})\nonumber \\
&=& \begin{pmatrix} 
{R_{1,3}\over R_{1,2} } &{R_{2,2}\over R_{1,2}R_{2,1} }& {R_{3,1}\over R_{2,1}R_{3,0} }& {1\over R_{3,0} }\\
0 &{R_{2,2}\over R_{2,1} }& {R_{1,2}R_{3,1}\over R_{2,1}R_{3,0} }&  {R_{1,2}\over R_{3,0} }\\
0 & 0 & {R_{3,1}\over R_{3,0} } & {R_{2,1}\over R_{3,0} }\\
0 & 0 & 0 & 1
\end{pmatrix}
 \label{matU}
\end{eqnarray}
\end{example}

Let us apply the result of Theorem \ref{solonear}, with $k_0=m_1+1=k_1$ and $k=n+m_1$, $j_0=j-n+1$,
$j_1=j+n-1$. Defining $\epsilon=1-(k\, {\rm mod}\, 2)$ and $\theta=j_1$ mod $2$, we get:
\begin{eqnarray*}
R_{1,n+m_1} &=&T_{1,\epsilon,k}=T_{1,\theta,k_1} \left(({\tilde U}_\bm U_\bm)^{n-1}\right)_{1,1}\\
&=& R_{1,m_1+1} \left(({\tilde U}_\bm U_\bm)^{n-1}\right)_{1,1}
=R_{1,m_1}\left( (U_\bm {\tilde U}_\bm)^n\right)_{1,1}
\end{eqnarray*}
where we have used $(U_\bm)_{1,x}=\delta_{1,x}$ and 
$({\tilde U}_\bm)_{x,1}=\delta_{x,1}T_{1,1-\theta,k_1-1}/T_{1,\theta,k_1}=\delta_{x,1}R_{1,m_1}/R_{1,m_1+1}$.
Comparing this with \eqref{tobecomp}, we deduce the

\begin{thm} \label{resoident}
Let $\bm$ be a Motzkin path of length $r-1$, $T_\bm$ the $2r+2\times 2r+2$ transfer matrix of the
$Q$-system solution of Ref. \cite{DFK3}, and $U_\bm$, ${\tilde U}_\bm$ as above. Then we have an identity
between ``resolvents":
$$ \left((I-T_\bm)^{-1}\right)_{1,1}=  \left((I-t U_\bm{\tilde U}_\bm)^{-1}\right)_{1,1} $$
\end{thm}

Note that the matrices $U_\bm$ and ${\tilde U}_\bm$ have size $r+1\times r+1$. 
We may however view the product
$U_\bm{\tilde U}_\bm$ as the transfer matrix of a weighted graph $G_\bm$ with
$2r+2$ vertices labeled $1,2,...,(r+1)$, $1',2',...,(r+1)'$, defined as follows.
We interpret the matrix element $(U_\bm)_{i,j}$ as coding the weights of the oriented edge $i\to j'$
of $G_\bm$
while $({\tilde U}_\bm)_{i,j}$ codes the weight of the oriented edge $i'\to j$. Alternatively, we may form the
$2\times 2$
block transfer matrix $\theta_\bm=\begin{pmatrix} 0 & {\tilde U}_\bm \\ U_\bm & 0\end{pmatrix}$ for $G_\bm$
and note that $\Big((I-\sqrt{t}\, \theta_\bm)^{-1}\Big)_{1,1}=\Big((I-t U_\bm{\tilde U}_\bm)^{-1}\Big)_{1,1}$.

Noting that the non-zero elements of 
${\tilde U}_\bm$ have the same indices as those of the transpose of $U_\bm$, we see that $G_\bm$
has doubly oriented edges only, with specific weights for each orientation. From the network construction,
all these weights are Laurent monomials of the initial data.

\begin{example}
For the Motzkin path $\bm=(2,1,0)$ of $r=3$, we may use the matrices \eqref{matU}. Without altering
the resolvent, we may gauge-transform the matrices $U_\bm$ and ${\tilde U}_\bm$ with invertible
matrices $R,L$ with $L_{j,1}=L_{1,j}=\delta_{j,1}$:
$V_\bm=L^{-1}U_\bm R$ and ${\tilde V}_\bm=R^{-1}{\tilde U}_\bm L$.
For the choices:
$$L={\rm diag}(1,{R_{1,2}R_{2,1}\over R_{2,2}},{R_{2,1}R_{3,0}\over R_{3,1}}, R_{3,0})
\qquad R={\rm diag}(1,R_{1,2},R_{2,1},R_{3,1})$$
we find that
$$V_\bm=\begin{pmatrix}
1 & 0 & 0 & 0 \\
y_2 & y_3 & 0 & 0\\
y_{3,1} & y_4 & y_5 & 0 \\
y_{4,1} & y_{4,2} & y_6 & 1
\end{pmatrix}\qquad 
{\tilde V}_\bm=\begin{pmatrix}
y_1 & 1 & 1 & 1\\
0 & 1 & 1 & 1\\
0 & 0 & 1 & 1 \\
0 & 0 & 0 & y_7
\end{pmatrix}$$
in terms of the parameters of \eqref{params}. 
One checks directly the statement of Theorem \ref{resoident} with the matrix \eqref{tmatq}
and by computing the rational
fraction $\Big( (I-tV_\bm {\tilde V}_\bm)^{-1}\Big)_{1,1}$.
\end{example}

Finally, turning to $R_{\al,n}$ for $\al>1$, we have a simpler picture than in \cite{DFK3},
as the quantity $R_{\al,n+m_1}/(R_{\al,m_1})^\al$ is expressed directly as
the partition function for $\al$ non-intersecting paths on the associated network, without having
to generalize the Lindstr\"om-Gessel-Viennot Theorem.

\section{Conclusion}

In this paper we have presented an explicit solution of the $A_r$ $T$-system in terms
of arbitrary boundary data. This solution is expressed in terms of partition functions
of weighted paths on some particular networks, determined by the boundary.

We have briefly described the connection of the $T$-system to a particular cluster algebra.
In particular, the sets of boundary data we have considered here form only a subset of the clusters
in this cluster algebra. What happens is that the form of the cluster mutations can change in
general from the equation \eqref{Tsys}, in some sense, the equation itself evolves, 
leading to clusters of a different kind.
Nevertheless, the positivity conjecture seems to hold for these other clusters as well.
It would be extremely interesting to probe whether these other clusters also have a description
in terms of networks, that would make the Laurent positivity property manifest, like in the cases
studied in this paper.

Another possible direction of generalization concerns non-commutative cluster algebras. 
In \cite{DFK09b}, the cluster algebra for the non-commutative $A_1$ $Q$-system was introduced,
and its positivity proved by use of a non-commutative weighted path model. We have checked that 
it can be reformulated as a non-commutative network model in the spirit of the present paper,
but with transfer matrices with entries in a non-commutative algebra. We hope to report on this 
direction in a later publication.

\medskip
\noindent{\bf Acknoledgments.} We would like to thank R. Kedem for numerous discussions.
We also thank S. Fomin for hospitality at the University of Michigan
and many discussions while this work was completed.
We received partial support from the ANR Grant GranMa, the
ENIGMA research training network MRTN-CT-2004-5652,
and the ESF program MISGAM.


\end{document}